\documentclass[preprint,3p,authoryear]{elsarticle}
\usepackage[utf8]{inputenc}
\usepackage{multirow}
\usepackage{booktabs}
\usepackage{eurosym} 
\usepackage{url}
\usepackage{graphicx}
\graphicspath{ {figures/} }
\usepackage{color}
\usepackage{amsmath}
\usepackage{pifont}
\newcommand{\cmark}{\ding{51}}
\newcommand{\xmark}{\ding{55}}
\usepackage[acronym,shortcuts]{glossaries}
 \newacronym{CBA}{CBA}{Cost–Benefit Analysis}
 \newacronym{CDNLP}{CDNLP}{Continuous Dynamic Network Loading Problem}
\newacronym{IM}{IM}{Infrastructure Manager}
\newacronym[firstplural = Railway Undertakings (RUs)]{RU}{RU}{Railway Undertaking}
\newacronym[firstplural = Freight Operating Companies (FOCs)]{FOC}{FOC}{Freight Operating Company}
\newacronym{TAC}{TAC}{Track Access Charge}
\newacronym{PSO}{PSO}{Particle Swarm Optimization}
\newacronym{GA}{GA}{Genetic Algorithms}
\newacronym{IP}{IP}{Interior-point Method}
\newacronym{AS}{AS}{Active-set Method}
\newacronym{SQP}{SQP}{Sequential Quadratic Programming}
\newacronym{PS}{PS}{Pattern Search}
\newacronym[firstplural = Regulatory Bodies (RBs)]{RB}{RB}{Regulatory Body}
\newacronym{NPV}{NPV}{Net Present Value}
 
\newacronym{WTP}{WTP}{Willingness-To-Pay}
\newacronym{WTT}{WTT}{Whole-arc Travel Time}
 
\usepackage[linesnumbered,ruled,vlined]{algorithm2e}
\usepackage{setspace}

\SetCommentSty{mycommfont}
\SetKwInput{KwInput}{Input}                
\SetKwInput{KwOutput}{Output}

\date{}

\begin{document}

\begin{frontmatter}

\title{A model for pricing freight rail transport access costs: economic and environmental perspectives
}

\author[inst1]{Ricardo García-Ródenas}
\author[inst2]{Esteve Codina} 
\author[inst3]{Luis Cadarso}
\author[inst1]{\\María Luz López-García}
\author[inst1]{José Ángel Martín-Baos\corref{cor1}}
\cortext[cor1]{Corresponding author}
\ead{JoseAngel.Martin@uclm.es}

\affiliation[inst1]{organization={Department of Mathematics, Escuela Superior de Inform\unexpanded{á}tica},
      addressline={University of Castilla-La Mancha}, 
      city={Ciudad Real},
      postcode={13071}, 
      country={Spain}}
\affiliation[inst2]{organization={Department of Statistics and Operations Research},
      addressline={Universitat Politècnica de Catalunya}, 
      city={Barcelona},
      postcode={13071}, 
      country={Spain}}
\affiliation[inst3]{organization={Aerospace Systems and Transport Research Group, European Institute for Aviation Training},
	addressline={Rey Juan Carlos University}, 
	city={Fuenlabrada, Madrid},
	postcode={28942}, 
	country={Spain}}
  
\begin{abstract}

In deregulated railway markets, efficient management of infrastructure charges is essential for sustaining railway systems. This study sets out a method for infrastructure managers to price access to railway infrastructure, focusing on freight transport in deregulated market contexts. The proposed methodology integrates negative externalities directly into the pricing structure in a novel way, balancing economic and environmental objectives. it develops a dynamic freight flow model to represent the railway system, using a logit model to capture the modal split between rail and road modes based on cost, thereby reflecting demand elasticity. The model is temporally discretized, resulting in a mesoscopic, discrete-event simulation framework, integrated into an optimization model that determines train path charges based on real-time capacity and demand. This approach aims both to maximize revenue for the infrastructure manager and to reduce the negative externalities of road transport. The methodology is demonstrated through a case study on the Mediterranean Rail Freight Corridor, showcasing the scale of access charges derived from the model. Results indicate that reducing track-access charges can yield substantial societal benefits by shifting freight demand to rail. This research provides a valuable framework for transport policy, suggesting that externality-sensitive infrastructure charges can promote more efficient and sustainable use of railway infrastructure.

\end{abstract}




\end{frontmatter}

\section{Introduction}

The European Directive EU 91/440/EEC \citep{LC1} specified guidelines for the deregulation of railway markets across Europe, advocating for a governance model characterized by vertical separation of the railway sector (refer to \cite{ElE22} for a detailed analysis). This separation involved creating two separate entities: the \gls{IM}, responsible for the provision and management of railway infrastructure, and \glspl{RU}, which utilize this infrastructure to offer freight or passenger services. \glspl{RU} submit requests for train paths to the \gls{IM}, who allocates these paths and collects the corresponding \glspl{TAC} from the network users, i.e., the RUs. Track access charges can be highly significant for rail operators, in some cases accounting for up to $88\%$ of their costs \citep{ScK22}.

The adoption and delineation of charging principles and the calculation of \glspl{TAC} remain key topics of discussion among \glspl{RB}, states, \gls{IM}, and \glspl{RU} (see, e.g., \cite{EU20}). Two main theoretical pricing principles are applied across European countries. The first, based on the ``institutional school'' of economics, asserts that fees should reflect the actual costs of producing and using services, without artificial price distortions from state intervention. The second, rooted in the ``marginalist school'', suggests that pricing should follow supply and demand principles, setting prices to reflect the marginal cost of providing an additional service, alongside consumers' willingness to pay for it. These approaches have significant implications not only for the profitability of \glspl{IM} and \glspl{RU}, but also for the overall efficiency and sustainability of the transport system.

European regulation has introduced legislation on infrastructure charges, beginning with Directive 2001/14/EC \citep{LC2}, which advocates for short-run marginal cost pricing. However, the regulations also allow \glspl{IM} to apply mark-ups for full infrastructure cost recovery, but only ``when the market can sustain it''.

The mark-up principles employed by most European \glspl{IM} often conflict with strategic decisions to build or expand railway infrastructure. Following a strictly profit-based approach focused on maximizing financial returns could lead to the conclusion that large-scale investments are infeasible, as they may not yield sufficient economic returns alone. This highlights the need for access charges criteria to align with societal goals, considering both monetary and non-monetary benefits to maximize total societal welfare.  This study is positioned within the framework of tactical planning and examines pricing mechanisms for railway infrastructure access, considering both the level of infrastructure utilization and the negative externalities associated with freight transport.
    
Furthermore, this issue is pivotal in meeting the European Union (EU) targets for shifting freight transport from road to more sustainable modes. The EU aims for a $30\%$ modal shift by 2030 and at least $50\%$ by 2050 for shipments exceeding 300 km. Rail freight transport is recognized as a cost-effective and environmentally sustainable alternative due to its economies of scale, reduced emissions, and ability to mitigate other negative externalities.

This study focuses exclusively on rail freight transport, with \glspl{RU} represented by \glspl{FOC}, and assumes that passenger rail services are given priority, while freight services are allocated based on residual capacity. The research question addressed is how to design \gls{TAC}  schemes that simultaneously incorporate cost recovery objectives, available infrastructure capacity, market conditions, and the full range of social benefits associated with rail infrastructure use. In this context, an optimization approach is proposed for pricing train paths (in terms of infrastructure charges), enabling the \gls{IM} to maximize access charge revenue while accounting for negative externalities such as environmental impact and road safety.

The proposed methodology presents several key challenges. Firstly, the model must account for demand elasticity to capture the effects of costs and travel times on modal split. To simplify this, road freight transport is represented through a logit-based modal split model. Secondly, instead of a conventional per-kilometer toll scheme, the approach implements capacity-based pricing, where charges vary based on network congestion levels, imposing higher fees on heavily used routes compared to corridors with lower demand. This is achieved through a dynamic freight flow model incorporating arc capacities, temporal demand variations, and capacity fluctuations (e.g., day vs. night operations).

This methodology has been implemented to ascertain \gls{TAC} for the Mediterranean Rail Freight Corridor (RFC6). In this context, two experiments are proposed. The first experiment assesses the proposed \gls{TAC} calculation framework, whereas the second experiment explores the broader discourse on whether \glspl{TAC} ought to be designed exclusively for cost recovery or as a strategic policy instrument for guiding railway infrastructure investment. Specifically, it examines capacity allocation and infrastructure pricing policies at the tactical level as essential components of long-term investment strategies aimed at fostering a sustainable transport system. This perspective regards \gls{TAC} not only as a financial mechanism but also as a means for attaining sustainability objectives in rail freight transport.

\section{Literature review and contributions}

Effective capacity allocation and infrastructure pricing policies are essential for ensuring both the efficiency and equity of railway network operations. Unlike other resources, rail capacity cannot be stored, making \gls{TAC} the primary mechanism for aligning network utilization with available capacity. Consequently, capacity allocation and pricing must be addressed jointly, as their interdependence adds significant complexity to the problem. Understanding how these policies shape service levels and network performance is particularly important in vertically separated rail markets, where infrastructure managers must balance competing demands from passenger and freight services.

In vertically separated rail markets, the \gls{IM} is responsible for capacity management, a task complicated by the heterogeneous nature of rail traffic. Passenger and freight services compete for network capacity while operating under distinct conditions, such as differences in speed, scheduling requirements, and flexibility of service. Passenger demand is highly schedule-sensitive, requiring carefully coordinated timetables, whereas freight services exhibit greater flexibility, adapting to variations in demand and operational constraints. Even within passenger transport, significant distinctions exist between commercial services with long-term planning horizons, subsidized public transport services, and short-term traffic, which requires greater operational adaptability.  

Similarly, the welfare effects of open-access competition in railway markets differ for passenger and freight sectors. Calculating societal benefits for passenger rail services involves assessing changes in commuter travel times, waiting times, transfers, crowding, and operating costs for commuter train operators. For freight, it is essential to consider the externalities of freight transport on the economy, environment, climate, and society. By 2005, the components considered for access fees across 23 European countries were documented in \cite{Nash2005}, with few including externalities such as accidents and environmental impacts. Since then, the inclusion of externalities in access fee structures has evolved. \cite{DHS15} identify several key negative externalities, including greenhouse gas emissions, water pollution, and land use, which primarily impose external costs on the environment and climate, as well as air pollution, noise, congestion, and accidents, which negatively affect human health and create additional economic burdens. Moreover, \cite{CMS21} and \cite{BGE19} compare the costs of various negative externalities across different freight transport modes for 2019. With a rail-to-road freight transport ratio of $1:5$ in tonnage, the associated negative externality costs show a ratio of $1:13$. These disparities strongly support policies promoting a modal shift from road to rail in the movement of freight. An effective access cost structure can further encourage rail demand by improving its competitiveness with respect to other transport modes.

\cite{Gib03,BEA22} identify three fundamental approaches to railway capacity allocation, which can be collectively referred to as capacity allocation mechanisms. These include: i) administrative methods, where capacity is assigned based on predefined rules and priority criteria; ii) cost-based approaches, which allocate capacity considering cost structures and externalities, often through social cost-benefit analysis; and iii) market-based mechanisms, where willingness-to-pay principles guide allocation through competitive bidding or pricing strategies. 

To resolve capacity conflicts, administrative methods are commonly applied, relying on straightforward, predetermined criteria such as prioritizing faster trains over slower ones or giving precedence to passenger trains over freight. These procedures are relatively simple to implement, quick, predictable, and transparent, making them suitable for monopoly environments. However, as \cite{HZM13} highlights, in transitioning to liberalized markets, new entrants often face barriers imposed by incumbents, infrastructure managers, rail regulators, and terminal operators. To foster effective competition in liberalized railway markets, \glspl{IM} must guarantee non-discriminatory infrastructure access through fair, transparent, and manageable allocation methods. While administrative mechanisms meet these requirements, they lack an inherent pricing mechanism for access charges and thus fall outside the scope of this study.

In market-based allocation systems, \glspl{RU}'s capacity requests are prioritized based on their \gls{WTP}. This approach maximizes societal value by ensuring that limited resources, such as congested railway tracks, are allocated efficiently through a pricing mechanism reflecting the economic valuations of users. Several versions of market-based allocation exist. A common method is static congestion charging, which reduces demand to align with available capacity by imposing higher charges during peak congestion periods. Another approach is dynamic pricing, which continuously adjusts track access charges based on factors such as booking lead time, current capacity utilization, and demand fluctuations. Dynamic pricing optimizes resource allocation by aligning prices with real-time conditions, effectively managing demand while incentivizing early bookings and encouraging operators to select less congested time slots. This approach enhances resource utilization, mitigates peak-time congestion, and fosters a more balanced and efficient transportation system.

From the perspective of the ``institutional school'', economic efficiency is achieved through resource allocation guided by the cost-benefit principle. This approach emphasizes that prices should reflect both private and social costs, including externalities such as environmental impacts, to ensure optimal resource allocation. Furthermore, it is argued that infrastructure charges should be equitably distributed across different transport modes, a dimension often neglected in earlier approaches.
\gls{CBA} is a widely-used framework in transportation science, often employed to evaluate potential infrastructure investments and transport policies, including pricing strategies. This method systematically quantifies all costs and benefits in monetary terms to assess net impacts. In the context of transport policy analysis, key components typically include consumer surplus (benefits accrued by passengers and freight customers), producer surplus, externalities, and net public expenditure. For instance, \cite{PrT18} apply \gls{CBA} not to assess railway investments but to evaluate capacity allocation and pricing strategies. This application is particularly noteworthy as it incorporates externalities as a core element within the structure of costs and benefits, highlighting their critical role in transport policy evaluation.

The train path allocation problem is addressed in \cite{Klabes} for the European context. They formulated the problem conceptually using the blocking time theory to represent capacity, stating the interactions with existing train timetable construction tools and proposed heuristically-based algorithmic solutions including as a criterion the dissatisfaction of \gls{RU}'s in a game-theoretical framework, without taking externalities into account.

\cite{Har13} gives a comprehensive classification of train-path pricing methods for freight transport, emphasizing the role of economic criteria in decision making by the \gls{IM}. These methods can be categorized as follows:
\begin{enumerate}
    \item[-] {\sl Fixed Allocation.} Train paths are priced based on predetermined accounting rules, such as distance traveled, or tonnage carried. This approach aims to distribute both fixed and marginal costs equitably among trains or routes, providing a straightforward and transparent allocation mechanism.
    \item[-] {\sl Marginal Costs.} Similar to fixed allocation schemes, but focusing solely on marginal costs. This method ensures that prices reflect the additional cost of providing a train path, encouraging cost-efficient usage of the infrastructure.
    \item[-] {\sl Value of Service.} Train paths are priced according to their estimated relative value of service. This valuation is typically influenced by the cargo's economic value, allowing higher-value goods to justify higher access charges.

    \item[-] {\sl Yield Management.} The \gls{IM} adopts a dynamic and confidential pricing strategy, often employing price discrimination based on \glspl{FOC}'  estimated \gls{WTP}. Prices are adjusted according to congestion levels and game-theoretic considerations, treating train paths as perishable commodities. This approach encourages network use during periods of low congestion and maximizes revenue during peak times.
    \item[-] {\sl Auction.} Train paths are allocated through an open auction or bidding process. In this method, the \gls{IM} exerts limited control over individual buyer prices, and \glspl{FOC} gain visibility into the prices paid by other participants. This competitive environment aims to achieve market-driven allocation and pricing of resources.
\end{enumerate}

The previous methods, based on mark-up principles aim to reduce government subsidies by transferring infrastructure costs to \glspl{FOC}. The first three methods in the list use specific cost rules to determine these charges, focusing on recovering investment or charging for the value of the service, without necessarily maximizing revenue. In contrast, ``Yield Management'' and ``Auction'' methods are designed to maximize aggregate revenue by setting prices according to operators' willingness to pay, capturing additional revenue beyond basic cost recovery. However, none of these approaches directly capture network service levels, although ``Yield Management'' and ``Auction'' do so in an indirect manner.



In an effort to address traffic heterogeneity, \cite{BEA22} propose a hybrid method that combines \gls{WTP} and \gls{CBA} for capacity allocation, allowing the approach to be tailored to different market segments. \cite{Stojadinovic} examine train path performance (in line with Directive 2012/34/EU; \cite{LC3}) and analyze capacity requests from various operators as an alternative to relying on previous timetables.

Table~\ref{tab:state-of-the-art} presents a structured review of key studies on railway access pricing and capacity allocation, facilitating a comparative analysis of their objectives, pricing methodologies, network characteristics, and consideration of externalities. The first three columns capture the reference, study objective (e.g., efficiency improvement, capacity allocation), and pricing strategy, distinguishing between willingness-to-pay, cost-benefit analysis, or hybrid approaches. The {\sl Railway Infrastructure} category is subdivided into three dimensions: {\sl Network Type}, specifying physical characteristics (e.g., single-track lines, corridors); {\sl Aggregation Level}, indicating the granularity of analysis (microscopic, mesoscopic, or macroscopic), which is particularly important for assessing network capacity and how railway conflicts are analyzed; and {\sl Modeling Approach}, describing the methodological framework. Finally, the {\sl Externalities} column identifies whether the study explicitly incorporates social and environmental costs into the pricing framework. This structured classification ensures a clear and systematic comparison of methodologies, highlighting how different approaches integrate economic and sustainability considerations into railway infrastructure pricing.

\begin{table}[h]
\centering
\caption{Summary of Key Studies on Capacity Allocation and Infrastructure Pricing. \label{tab:state-of-the-art}}
\resizebox{1\textwidth}{!}{
\begin{tabular}{p{3cm} p{3cm} p{3cm} p{2.5cm} p{2.5cm} p{4cm} c }
\toprule
{\bf Reference} & {\bf Objective} & {\bf Capacity Allocation} & \multicolumn{3}{c}{\bf Railway Infrastructure} & {\bf Externalities} \\
\cmidrule(lr){4-6}
 &  &  & {\bf Network Type} & {\bf Aggregation Level} & {\bf Modeling Approach} &  \\
\midrule
\cite{Nilsson} & Efficient capacity allocation and social welfare enhancement & \gls{WTP} (Auction) & Single-track line & Macroscopic & Dual optimization with auction-based scheduling & \xmark\\
\cite{Klabes} & Railway capacity allocation under competitive markets & \gls{WTP} (Game-theoretical approach) & Corridors & Macroscopic & Nash equilibrium-based allocation with congestion assessment & \xmark\\
\cite{Har13} & Transparent and revenue-maximizing train path allocation & \gls{WTP} (Auction) & Single-track line & Microscopic & Iterative auction-based scheduling with bid round optimization & \xmark\\
\cite{BGL06} & Efficient capacity allocation through auction-based scheduling & \gls{WTP} (Combinatorial Auction) & Double-track lines & Macroscopic & Iterative combinatorial auction with integer programming & \xmark\\
\cite{BuB14} & Universal access charge model based on cost recovery & Systems approach & Mid- and small-size railway networks & Macroscopic & Cost-driven segmentation and modular tariff optimization & \xmark\\
\cite{LLC14} & Capacity charging for shared-use rail corridors & \gls{CBA} & Freight and passenger corridors & Macroscopic & Congestion and opportunity cost estimation & \xmark\\
\cite{Ali2020} & Societal cost-based pricing of commercial train paths & \gls{WTP} (Yield Management) & Mixed traffic (commuter and commercial trains) & Microscopic & Societal cost assessment and reservation price calculation & \xmark\\

\cite{MTP18} & Incentive-based access charge reduction for freight operators & \gls{WTP} & Single line & Microscopic & Infrastructure gap compensation model & \xmark\\

\cite{PrT18} & Societal impact evaluation of capacity allocation and pricing policies & \gls{CBA} & Freight and passenger corridors & Macroscopic & Scenario-based cost-benefit analysis & \cmark\\

\cite{BRL21} & Welfare-optimal track access charges for (passenger) rail & Welfare optimization & Single segment & Aggregated demand model (Nested Logit Model) & Pricing and frequency optimization & \cmark\\

\cite{BEA22} & Integrated capacity allocation for heterogeneous traffic & \gls{WTP}+\gls{CBA} & Predefined train paths & Macroscopic & Auction-based allocation, social CBA, and dynamic pricing & \cmark\\

This paper & Externality-sensitive access charge optimization for freight rail & \gls{WTP}+\gls{CBA} & Double-track line & Mesoscopic & Dynamic freight flow simulation with modal split modeling & \cmark\\
\bottomrule
\end{tabular}
}
\end{table}

\subsection{Main contributions}

This study examines the complex issue of track access charges in liberalized railway markets by developing a pricing model that integrates a dynamic traffic loading approach with capacity constraints, a modal-split model, and environmental considerations. The train path allocation is framed within an operational planning context, where rail network conflicts must be explicitly accounted for to accurately determine available capacity. Given that track access charges are addressed at a tactical level in this study, capacity considerations are essential. Consequently, capacity has been modeled from a mesoscale perspective using a dynamic traffic loading model, which is appropriate for the tactical nature of pricing decisions and explicitly accounts for network capacity.

The proposed framework leads to an optimal control problem, where discretization transforms it into an optimization model with constraints defined by a discrete event simulation, parameterized by track access charges. Due to the complexity of this problem, various optimization algorithms have been tested using the Mediterranean Rail Freight Corridor (RFC6) as a case study. These include gradient-based methods (Sequential Quadratic Programming, Interior Point, and Active Set methods) and function evaluation-based approaches (Genetic Algorithms, Particle Swarm Optimization, and Pattern Search). Pattern Search emerges as the most efficient, although it remains computationally demanding.

The main contributions of this study to the literature are as follows:
\begin{enumerate}
	\item {\sl Integration of modal competition in \gls{TAC} design.} This study proposes a novel methodology for designing \glspl{TAC} that explicitly accounts for their impact on modal competition with road transport. By incorporating this effect, the model provides a more accurate assessment of transport externalities, addressing a critical gap in existing pricing frameworks.
	\item {\sl Holistic optimization of \glspl{TAC} and train-path interdependencies.} The model accounts for the non-additive nature of infrastructure costs and the interdependencies between train paths sharing railway segments. Rather than treating network segments independently, it optimizes total revenue from train path allocation over the entire planning horizon, ensuring a more efficient pricing strategy. Additionally, the model is adaptable to a conventional per-kilometer toll scheme, offering a clearer and more transparent pricing structure for operators. Numerical experiments show that this scheme yields comparable outcomes to the holistic approach, highlighting its practicality and ease of implementation.
	\item {\sl Economic and environmental trade-offs in \gls{TAC} design.} Numerical experiments on the Mediterranean Rail Freight Corridor underscore the importance of externality-based \glspl{TAC} in balancing economic returns for infrastructure managers with environmental benefits. Furthermore, the study highlights the dual role of \glspl{TAC} as both cost recovery mechanisms and investment tools, emphasizing the need to consider broader societal impacts when setting tariffs. One  contribution of this work is the development of a tool to facilitate this analysis, which is made openly available to researchers. This approach aligns with EU objectives for sustainable and energy-efficient transportation systems, offering a valuable resource for future studies in this area.
\end{enumerate}

\section{Optimization model for setting track access charges}

The model must effectively capture two key attributes of the rail network: the capacities of track segments and the dynamic nature of flows and  travel times. Firstly, this paper adapts the \gls{CDNLP}, \cite{XWF99}, originally developed for traffic assignment problems, to address these specific requirements. The \gls{CDNLP} computes link and path travel times within a congested network by using time-dependent path flow rates over a specified period.  Secondly, similar to \cite{RoC22}, the proposed model features,  in a dynamic context, a logit-based demand model for the modal split of the freight between rail and road, depending on transport costs. These costs are influenced by the \gls{TAC} and travel times, thereby shaping the operational expenses of the \glspl{FOC}. Finally, the objective function includes the \gls{TAC} criterion used by the \gls{IM}, serving as a proxy for revenue management. The optimization variables are the access charges for each train path at each time. Each of these elements is described in detail below.

Table~\ref{tab:notation} provides a comprehensive list of notations used to describe the elements involved in path flows, path composition, and model demand. 


\begin{table}[!]
\centering
\begin{tabular}{l@{\hspace{1.em}
    \hspace{1em}}p{0.75\linewidth}}
	\toprule 
	\multicolumn{1}{l}{\textbf{Indices and sets}} &
	\multicolumn{1}{@{}l}{\textbf{}}\\%
	\midrule 
	$a:$ & Arc in the railway network.\\
	$r:$ & Path on the railway network.\\
	${a^+_r}:$ & Arc immediately following arc $a$ on  path $r$.\\
	${a^-_r}:$ & Arc immediately preceding arc $a$ on  path $r$.\\
	$o_r:$     &  First arc (connector) of  path $r$.\\
	$d_r:$     & Last arc (connector) of  path $r$.\\
    $W:$       & The set of active origin-destination pairs on the network.\\ 
    $V:$       & Set of vertices in the railway network. \\ 
    $A:$       & Set of arcs in the railway network. \\  
    $A_r:$     & Set of arcs that constitute path $r$. \\  
    $R:$       & Set of  paths in the railway network. \\  
    $R_\omega:$ & Set of  paths corresponding to the o-d pair $\omega$. \\  
    $R_a:$      & Set of paths that use arc $a$. \\  
    \rule{0pt}{1em} & \\ 
    \midrule
	\multicolumn{1}{l}{\textbf{Variables}} &
	\multicolumn{1}{@{}l}{\textbf{}}\\%
	\midrule
    $t:$ & Time (independent variable)\\
	$f_a^\ast (t):$ & Entry flow intensity for the arc $a$ at time $t$ (measured in trains per unit of time).\\
	$f_{ar}^\ast (t):$ & Entry flow intensity for the arc $a$ associated with the path $r$.\\
	$\phi_a (t):$    & Flow intensity exiting the  queue at arc $a $.\\
    $\phi_{ar} (t):$  & Flow intensity exiting the  queue at arc $a$ associated with the path $r$.\\
	$f_a (t):$    & Flow intensity served at the end of arc $a$.\\
	$f_{ar} (t):$ & Flow intensity served at the end of arc $a$ associated with the path $r$.\\
	$q_a(t):$     &  Vertical queue on arc $a$ representing unserved traffic (measured in train units). \\  
	$q_{r}(t)$:   & Number of vehicles on path $r$ at time $t$ (measured in train units). \\  
	$q_{ar}(t)$:  & Vertical queue on the arc  $a$ induced by path $r$. \\
	$\tau_r(t)$:  & Travel time to traverse  path $r$ departing at time $t$.\\
	$\nu_r(t)$:   & Time required for the last train to complete its journey on path $r$ at time $t$. \\  

	$D_\omega(t):$ & Demand intensity for pair $\omega$ at time $t$ (measured in tons per unit of time). \\  
	$D_{r}(t):$    & Demand intensity for path $r$ at time $t$ (measured in tons per unit of time). \\ 
    $L_\omega(t):$ & Demand intensity served by road for pair $\omega$ at time $t$ (measured in tons per unit of time). \\  
	$U_r(t):$ & The cost for the rail mode to transport one tonne-kilometer of freight using path $r$, starting at time $t$. \\  
	$\lambda_r(t):$ & Access charges for path $r$ at time $t$ per unit of flow (measured in monetary units per unit of time and vehicle flow). These are decision variables. \\
    \rule{0pt}{1em} & \\ 
    \midrule
	\multicolumn{1}{l}{\textbf{Parameters}} &
	\multicolumn{1}{@{}l}{\textbf{}}\\%
	\midrule
	$\delta_{ar}$:& Element of the  arc-path incidence matrix. It is  $1$ if arc $a$ is in path $r$, and $0$ otherwise.\\
	$\tau_a$:     &  Travel time for running section of the arc  $a$. \\
	$\tau^0_r$:   & Travel time to traverse train path $r$ with an empty network.\\
	$k_a(t)$:     & Arc capacity at instant $t$.\\
	$V_\omega:$   & The cost for road mode to transport one tonne-kilometer of freight for the origin-destination pair $\omega$.\\
	$\kappa:$    &  Conversion factor from freight flow to train flow for a train prototype in the network.\\
	$[0,T_{\max}]:$  & Planning horizon, i.e $t\in [0,T_{\max}]$. \\
	\bottomrule
\end{tabular}
\caption{Description of the notation used in the model}
\label{tab:notation}
\end{table}

\subsection{Dynamic network loading model}
The freight railway network is represented as a graph $G = (V, A)$, where $V$ is the set of vertices and $A$ is the set of edges. Each edge $a \in A$ represents a segment of rail infrastructure with homogeneous characteristics, such as loading gauge, maximum slope, permissible train weight, and maximum admissible train length. The dynamic network model captures the propagation of flows within the network and consists of two main elements: the basic arc model and the flow propagation model. The first element utilizes a point-queuing model similar to those applied in dynamic traffic assignment problems (e.g., \cite{LFK00,Garcia-Rodenas2006}). For simplicity, it is assumed that the railway network is empty at the initial time $t = 0$, and the operational period under consideration is $[0, T_{\max}]$. The planning horizon $T_{\max}$ is chosen to be sufficiently large so that the final state of the system is negligible relative to the total elapsed period.

\subsubsection{Arc-based model}

A crucial element of the \gls{CDNLP} is modeling arc performance. Two primary modeling frameworks are mainly used in arc loading models for traffic: the so-called \gls{WTT} models and the kinematic wave model of Lighthill–Whitham–Richards (LWR) for traffic flow. In this study, a \gls{WTT} model has been explored, which facilitates the development of a mesoscopic simulation model for railway networks. Conversely, LWR models are typically more suited to providing a microscopic characterization of railway networks.

\begin{figure}[!]
    \centerline{\includegraphics[scale=.4]{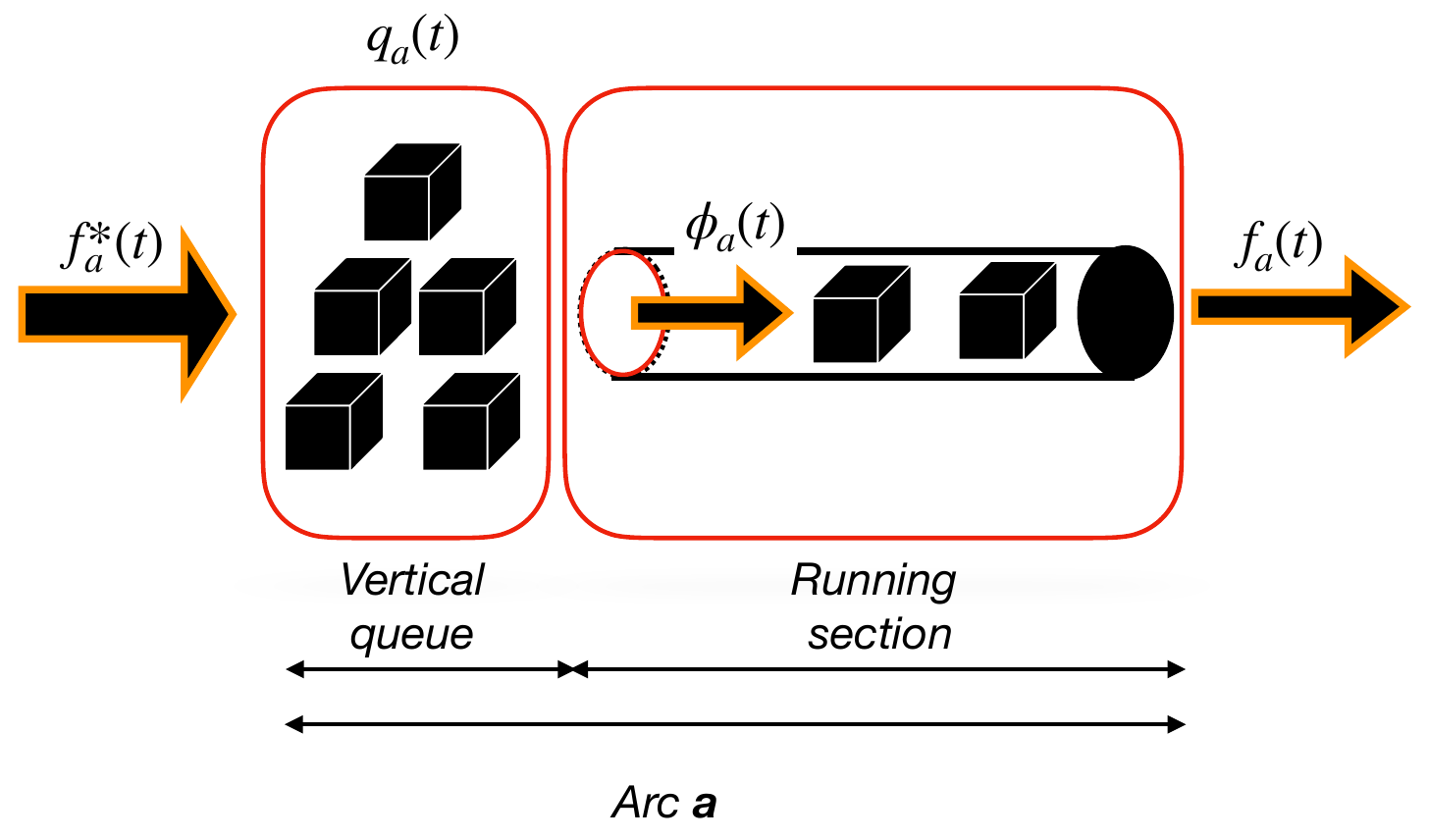}}
\caption{Representation of an arc $a$ at time instant $t$\label{fig:arco}}
\end{figure}

The \gls{WTT} queuing model employed in this work is illustrated in Figure~\ref{fig:arco} for a given arc $a$ at time instant $t$. It incorporates a vertical queue at the arc entrance, denoted by $q_a(t)$. The variables $f_a^*(t)$ and $\phi_a(t)$ denote the incoming and outgoing flows to/from the vertical queue, respectively, while $f_a(t)$ represents the outgoing flow from arc $a$. The travel time through the arc, denoted by $\tau_a$, is assumed to be independent of flow conditions, with congestion effects captured by the queue at the entry point. The network is modeled as a double-track system, where each track segment accommodates a unidirectional flow. Each arc is indexed by $a$.

The dynamics of flows on each arc $a \in A$ are governed by the following equations:
\begin{align}
    \label{eq1}
    q_a(t)&= \int_0^t  \left [f_a^\ast (\xi)-\phi_a(\xi)  \right ]{\rm d\xi}+q_a(0), & \forall a \in A, \\
    \label{eq2}
    \phi_a(t)&= \left \{
    \begin{array}{cl}
    \min \left \{f_a^*\left (t \right ),k_a \left (t \right ) \right \} \,\,\, &\hbox{ if }  q_a(t)=0\\
    k_a \left (t \right )   \,\,\, &\hbox{ if }  q_a \left (t \right )>0
    \end{array}
    \right .,  & \forall a \in A,\\
    \label{eq3:bis}
    f_a(t)&=\phi_a(t-\tau_a), & \forall a \in A. 
\end{align}

Equations (\ref{eq1}) and (\ref{eq2}) define the dynamics of flow and vertical queuing at the entrance of each arc $a$. Specifically, Equation~(\ref{eq1}) represents the evolution of the vertical queue, $q_a(t)$, based on the difference between incoming flow $f_a^\ast(t)$ and outgoing flow $\phi_a(t)$. Usually, the initial queue length, $q_a(0)$, is assumed to be zero. Equation (\ref{eq2}) restricts the outgoing flow $\phi_a(t)$ from the vertical queue to the arc capacity $k_a$ whenever a queue is present. If no queue is present, the outgoing flow equals the incoming flow, provided the arc capacity allows. Equation~(\ref{eq3:bis}) models the propagation of flow through the running section, assuming a congestion-independent travel time, $\tau_a$, from the start to the end of the arc.


\subsubsection{Network-based model}

The railway network is modeled as a directed graph consisting of arcs representing railway segments with homogeneous characteristics. These segments connect key locations, such as stations, depots, junctions, and modal interchange terminals, represented as nodes within the graph, referred to as regular nodes. Demand flows follow an origin-destination pattern, starting and ending at special artificial nodes known as \textit{centroids}. The set of active $o$-$d$ pairs is denoted by $W$, with an $o$-$d$ pair represented either as $(o, d)$ or simply by $\omega \in W$, depending on context. The set of paths connecting origin $o$ to destination $d$ is denoted by $R_\omega$, where $\omega = (o,d)$. 

Centroids are connected to other network nodes via special arcs called \textit{connectors}, which differ from the arcs representing physical railway segments—referred to as \textit{regular arcs}. Flows on the network arcs are described in terms of their path-flow composition, according to the paths linking the various origin-destination pairs.

The flows on network arcs are categorized into incoming flows, $f^a(t)$, and outgoing flows, $f_a(t)$, defined according to the network structure. Incoming flows $f^a(t)$ on an arc $a \in A$ may originate from external demand (through connectors) or from other adjacent arcs $a'$. Outgoing flows $f_a(t)$ represent flows exiting arc $a$, potentially serving as inputs to subsequent arcs in the network or as flows leaving the system. The relationship between arc flows and network paths is expressed using the binary parameter $\delta_{ar}$, indicating whether link $a$ is part of path $r$ ($\delta_{ar} = 1$) or not ($\delta_{ar} = 0$).

The flow dynamics on network arcs are primarily described by two types of equations: {\sl conservation equations} and {\sl propagation equations}. Conservation equations refer to flows $f^*_a(t)$ and flows $f_a(t)$ on network arcs $a \in A$ with a relationship due to the network structure. Flows $f^*_a(t)$ may consist of external inputs, either directly related to the demand or from other incident arcs $a'$ to arc $a$. Flows $f_a(t)$ are exiting flows from arc $a$ that will either reach its final destination or that will enter to vertical queue of emerging links $a'$ from $a$. The flows $f_a(t)$ and $f^*_a(t)$ are expressed in terms of the flows on the various paths $r$ that use arc $a$. The parameter $\delta_{ar}$ indicates whether path $r$ contains link $a$ ($\delta_{ar} = 1$) or not ($\delta_{ar} = 0$).

The flow conservation equations are as follows:
\begin{align}
    \label{eq3}
     f^*_{ar}(t) &=  f_{{a^{\hspace{-0.05cm}-}_r} r}(t), & \forall  a \in A_r, \; r \in R_\omega, \; \omega \in W,\\
     \label{eq4}
     f^*_a(t)&=\sum_{r \in R } \delta_{ar} f^*_{ar}(t), & \forall a \in A,\\
     \label{eq5}
     q_{ar}(t)&= \int_0^t \left [ f^\ast_{ar}(\xi)-\phi_{ar}(\xi) \right ] {\rm d\xi}, &  \forall a \in r, \; r \in R_\omega, \; \omega \in W,\\
     \label{eq6}
     q_a(t)&=\sum_ {r\in R} \delta_{ar} q_{ar}(t), & \forall a \in A.
\end{align}

Equation~(\ref{eq3}) specifies that the outflow from an arc on a given path $r$ matches the inflow on the subsequent arc. Next, Equation~(\ref{eq4}) states that the inflow in an arc $a$ is the sum of all the path inflows that traverse it. Then, Equation~(\ref{eq5}) establishes the queue in arc $a$ of trains traveling on path $r$. And finally, Equation~(\ref{eq6}) sums all queued trains for all paths that use the arc $a$ to generate the queue $q_a(t)$ at the arc.
    
The propagation of path flows on links $a\in A$, based on their path composition, is given by the following equations:
\begin{align}
\label{eq7}
\phi_{ar}(t)&=    \left \{
\begin{array}{ll}
\dfrac{f^*_{ar}(t)}{\sum_{r'\in R_a} f^*_{ar'}(t)} \cdot  \min \left \{f^*_{a}(t),k_a(t) \right \} \,\, &\hbox{ if } q_a(t) =0\\[0.15in]
\dfrac{q_{ar}(t)}{\sum_{r'\in R_a} q_{ar'}(t)} \cdot k_a(t), \,\, &\hbox{ if } q_a(t) >0
\end{array}
\right ., & \forall a \in A_r, \; r \in R_\omega, \; \omega \in W,\\
\label{eq8}
q_{r}(t)&= \int_0^{t} \left [ f^{\ast}_{o_r r}(\xi) -f_{d_r r}(\xi)\right ] {\rm d\xi}, & \; \forall  r \in R_\omega, \; \omega \in W,\\
\label{eq8:bis}
f_{ar}(t)&= \phi_{ar} (t-\tau_a), & \forall a \in A_r, \; r \in R_\omega, \; \omega \in W,\\
\label{eq9}
q_{r}(t)&= \int_t^{t+\nu_r(t)} f_{d_r r}(\xi){\rm d\xi}, & \; \forall r \in R_\omega, \; \omega \in W,\\
\label{eq10}
\tau_r(t)&=\max \left \{\tau^0_r,\nu_r(t) \right \}, & \; \forall r \in R_\omega, \; \omega \in W.
\end{align}

Equation~(\ref{eq7})  describes the dynamics of the vertical queue with path-disaggregated flows, i.e., the transformation of the Equation~(\ref{eq2}) given for homogeneous flows to path flows. It is assumed that the inflow of each path into the running section of an arc is directly proportional to the magnitude of its queue. Although this assumption deviates from the First-In-First-Out (FIFO) principle, the FIFO condition can be enforced as a criterion for queue management during the simulation of the model. Equation~(\ref{eq8:bis}) represents the propagation of disaggregated flows.

Equation~(\ref{eq8}) defines the queue for path flows, while Equations~(\ref{eq9}) and (\ref{eq10}) define $\tau_r(t)$ as the time required to clear all trains present on path $r$ at time $t$ or, equivalently, as the time needed for a train entering path $r$ at time $t$ to exit it. The rationale is that, assuming no overtaking, if a train catches up with its predecessor on path $r$, both will leave the system simultaneously. Therefore, the travel time on path $r$, $\tau_r(t)$, represents the time taken by a train to reach the end of path $r$ from its starting point and is always at least as long as the travel time when the path is empty. Note that when referring to the queue on a path, this can be at any point along the path, whether the train is in a running section or in a vertical queue.

\subsection{Demand model}

This section introduces a logit model designed to allocate freight transport between rail and road modes based on the associated costs for each mode. The network model considers various train types and freight categories. However, to maintain computational tractability, the model uses a single representative \textit{prototype} train and a single freight category. This simplification circumvents the necessity of introducing two additional indices into the formulation, one indexing train types and another freight categories.

The infrastructure performance for a given path $r$ is determined by the segment along that path with the most restrictive characteristics. This limitation affects train composition by setting a maximum allowable train length, consequently restricting both the number of freight cars and the total weight manageable by a single locomotive on path $r$. Therefore, once the prototype train is defined, the set of feasible paths $r$ operationally viable within the network is also determined.

The cost structure adopted is categorized as a service cost rather than a product cost, emphasizing that each transportation service incurs costs proportional to its length of use. The cost components are defined as follows:
i) $A_r(t)$, representing time-dependent costs such as those associated with train drivers, locomotives, and rolling stock;
ii) $\lambda_r(t)$, corresponding to payments of \glspl{TAC} made to the \gls{IM}, reflecting the infrastructure access charges imposed on railway operators. These charges significantly influence both the economic feasibility of rail transport and its competitive position relative to alternative transport modes;
and iii) $C_r$, denoting a fixed cost independent of train characteristics, covering expenses such as shunting operations at the origin or destination and other fixed administrative costs, expressed in units of ton-km.

The equations that define the modal split are as follows:
\begin{align}
\label{demand:eq1}
D_\omega(t)&=L_\omega(t) + \sum_{r\in R_\omega} D_r(t),& \forall \omega \in W,\\
D_{r}(t)&
= \dfrac{\exp \left (U_{r}(t) \right ) }{\exp(V_\omega)+\sum_{r' \in R_\omega}\exp \left (U_{r'}(t) \right )} \cdot  D_\omega(t), & \forall  r \in R_\omega,\, \omega \in W,
\label{demand:eq4}\\
\label{demand:eq3}
U_r(t)&= \left ( A_r(t) +  \lambda_r(t) \right ) \cdot \tau_r(t) +C_r, & \forall r \in R_\omega,\, \omega \in W, \\
\label{demand:eq2}
f^*_{o_r r}(t )&=  \kappa \cdot D_{r}(t), & \forall r \in R_\omega,\, \omega \in W.
\end{align}


Equation~(\ref{demand:eq1}) expresses that the total demand flow $D_\omega(t)$ is split into demand served by road, $L_\omega(t)$, and demand served by rail. The modal split is determined by a logit model, as shown in Equation~(\ref{demand:eq4}). The cost structures for rail transport are outlined in Equation~(\ref{demand:eq3}). Finally, Equation~(\ref{demand:eq2}) converts freight flows into train flows at the first arc $o_r$ of path $r$, enabling the flows to enter the network.

\subsection{Train Path Allocation Criterion}

The primary objective of the model is to determine optimal \glspl{TAC} within a tactical planning horizon. Consequently, the allocation of train paths is addressed as a secondary consideration, managed in a manner that effectively supports the primary goal.

Given the inherent interdependence between \glspl{TAC} and train path allocation, these two aspects ideally should be determined simultaneously. However, to prioritize pricing decisions, train path allocation is established indirectly. Specifically, it is assumed that passenger train services have predetermined schedules, thus establishing a residual capacity $k_a(t)$ available for freight transport at each time period $t$. Freight trains are dispatched only when demand meets a sufficient threshold to warrant their operation. Train paths are then allocated according to a criterion similar to stochastic user equilibrium assignment \citep{Pat94}, ensuring a probabilistic distribution of freight flows guided by cost minimization.

A fully microscopic model would be required to resolve potential conflicts at the train path level explicitly. Instead, the model incorporates these constraints in an aggregated manner by the use of vertical queuing on network arcs. If capacity constraints prevent immediate dispatch, freight trains may experience delays or be temporarily held at specified depots until movement is possible. It should be emphasized that the model does not differentiate between different \glspl{FOC}, as the principal focus remains on optimizing access charges rather than addressing specific operational constraints of individual operators.

\subsection{Objective Function}

The problem of determining \gls{TAC} for the railway network is inherently multi-objective. On one hand, the \gls{IM} aims to maximize revenue derived from access charges; on the other hand, as a publicly funded entity, it must also consider the negative externalities associated with freight transportation. The principal externalities in transportation include emissions (air pollution and greenhouse gases), noise pollution, water contamination, congestion, and accidents \citep{DHS15}.

The proposed objective function combines both goals, expressed as follows:
\begin{align}
\label{eq:objetivo}
    \hbox{Maximize }Z &=\sum_{\omega\in W}
    \sum_{r\in R_\omega }  \int_{0}^{T_{\max} }  \lambda_r(\xi) \cdot  \tau_r(\xi) \cdot  f_{o_r r}(\xi) \, {\rm d \xi} -
    \sum_{\omega \in W } \eta _\omega \int_{0}^{T_{\max} }  L_\omega (\xi) \,  {\rm d \xi} .
\end{align}

The first term of Equation~(\ref{eq:objetivo}) represents the total revenue received by the \gls{IM} from the use of train paths, reflecting dynamic pricing where different rates apply across the planning horizon $[0, T_{\max}]$. The second term accounts for transportation externalities, with the negative sign indicating the associated societal costs.

The formulation in the objective function relies on the assumption that externalities, denoted by $\varepsilon$, are dependent solely on the origin-destination pair $\omega$ and the total freight volume transported. Under this assumption, the externalities can be represented as follows:
\begin{align*}
    - {\cal E}&= - \left ({\cal E}_{truck}+{\cal E}_{train} \right )= - \left ( \sum_{\omega \in W} \eta^{truck} _\omega \int_{0}^{T_{\max} }  L_\omega (\xi) \, {\rm d \xi} 
 + \eta^{train} _\omega \int_{0}^{T_{\max} } [D_\omega (\xi ) -L_\omega (\xi) ] \, {\rm d \xi} \right )\\
 &= - \left (\sum_{\omega\in W} \underset{\eta_\omega}{\underbrace{(\eta^{truck} _\omega - \eta^{train}_\omega)}}\int_{0}^{T_{\max} }  L_\omega (\xi) \, {\rm d \xi} \right )+ ctant.
\end{align*}

In this formulation, the integral term quantifies freight volumes transported by road, and the parameter $\eta_\omega$ captures the specific characteristics of each origin-destination pair $\omega$, such as the distance covered by rail and road. This parameter quantifies the negative externalities of transporting a given amount of freight by road, compared to the externalities arising from transporting the same goods by rail.

To better understand the first term of the objective function, consider the case of constant pricing \(\lambda_r(t) = \lambda_r\) over the operating period \([0, T_{\max}]\), with a discretized timeline \(0 = t_0 < t_1 < t_2 \le t_{n_r} = T_{\max}\). Under the assumption that each discretized interval represents exactly one unit of train flow\footnote{It is assumed that the number of trains is an integer}:  
\[
\int_{t_{j-1}}^{t_j} f_{o_r r}(\xi) {\rm d} \xi = 1,
\]  
the following expression for \(\lambda_r\) will hold if the travel time \(\tau_r(t)\) is continuous (although this may not always be the case). Applying the mean value theorem for integrals, it is found that the total revenue received by the \gls{IM} from all trains operating on path \(r\), denoted by \(\Lambda_r\), can be calculated as:
\begin{align*}
    \Lambda_r=& \sum_{j=1}^{n_r} \int_{t_{j-1}}^{t_{j}}  \lambda_r(\xi) \cdot  \tau_r(\xi) \cdot  f_{o_r r}(\xi) \, {\rm d \xi}=\sum_{j=1}^{n_r} \lambda_ r \int_{t_{j-1}}^{t_{j}}    \tau_r(\xi) \cdot  f_{o_r r}(\xi) \, {\rm d \xi} =\\
    & = \sum_{j=1}^{n_r} \lambda_ r  \tau_r(c_j) \int_{t_{j-1}}^{t_{j}}   f_{o_r r}(\xi) \, {\rm d \xi}
    =\sum_{j=1}^{n_r} \lambda_ r  \tau_r(c_j)  \Rightarrow \lambda_r=\dfrac{\Lambda_r}{\sum_{j=1}^{n_r} \tau_r(c_j)}.
\end{align*}

The previous formula represents that the pricing intensity equals the total revenue divided by the total time of path usage. The previous expression can be rewritten as follows:
\begin{equation*}
    \lambda_r=\dfrac{ \Lambda_r}{\overline{\tau}_r\cdot n_r},
\end{equation*}
\noindent where \(\overline{\tau}_r\) represents the average travel time of the trains on path \(r\), and \(n_r\) is the number of trains using the  path $r$. Thus, the parameter \(\lambda_r\) can be interpreted as the income generated by a train per unit of time on path \(r\).

The model for setting \glspl{TAC} is formulated as an optimal control problem with the objective function given by Equation~(\ref{eq:objetivo}). Here, the time-dependent functions \(\lambda_r(t)\) act as decision variables (controls). The dynamics of the system are defined by Equations~(\ref{eq3})–(\ref{eq6}), with flow propagation constraints given by Equations~(\ref{eq7})–(\ref{eq10}) and modal split constraints in Equations~(\ref{demand:eq1})–(\ref{demand:eq2}).

\section{Resolution approach}

The resolution of the optimal control model formulated in the previous section is based on a discretization approach. Specifically, the model could be approached using two alternative methodologies for discretizing the freight flow model: time-based discretization and flow unit-based discretization. In the first approach, the planning horizon is divided into intervals \([t, t + \Delta t]\), within which freight flows are aggregated at each time step.  However, this method may lead to heterogeneous flow units, as the volume of freight transported within each interval can vary significantly.  To overcome this limitation, a flow unit-based discretization approach is taken, in which the discretization intervals are dynamically adjusted to ensure discrete and constant freight volumes equal to $\Delta f$. Each discretized flow unit is associated with a predefined quantity of freight and is represented by a conceptual entity referred to as a \textit{packet}, or simply \textit{train}. This approach not only maintains the required granularity for optimization but also provides a structured framework for defining movement rules within the network, as the discretization scheme determines how freight propagates across links and nodes.  

The discretization process extends to the objective function, where the continuous pricing functions \(\lambda_r(t)\) are transformed into a finite set of decision variables \(\lambda_{rj} = \lambda_r(t_j)\), corresponding to specific instants of time \(t_j\). Furthermore, the model constraints are reformulated within a discrete-event simulation framework, in which system dynamics evolve based on event-driven interactions. These events primarily correspond to packet arrivals at network nodes, either from incident arcs representing internal flow propagation or from external connectors modeling exogenous demand.  

The resulting simulation model is parameterized by the decision variables and serves as the fundamental mechanism for evaluating the objective function. Consequently, optimizing access charges requires integrating this simulation-based evaluation within an appropriate optimization algorithm. The following sections provide a detailed discussion of the discretization of the objective function and describe the simulation algorithm used to model system dynamics.

\subsection{Objective function}

It is assumed that the function $\lambda_r(t)$ can be approximated by a discretized piecewise-constant function over a grid of points $0= \delta_0< \delta_1< \ldots < \delta_k = T_{\max}$:
\begin{equation}
	\lambda_r(t):= 
	\begin{array}{cc}
		\lambda_{rj}& \text{ if } t\in [\delta_{j-1}, \delta_j).
	\end{array}
\end{equation}

The operational period $[0, T_{\max}]$ is discretized on a separate grid $0 = t_0 < t_1 < t_2 <\cdots < t_{n_r} \le  T_{\max}$, ensuring that the train intensity satisfies: 
\[
\int_{t_{j-1}}^{t_j} f_{o_r r}(\xi) \, {\rm d} \xi = \Delta f,
\]
which represents the packet size. 

The first term of the objective function~(\ref{eq:objetivo}), namely $\sum_{\omega} \sum_{r \in R_\omega}  \Lambda_r$, can thus be expressed as a function of the variables $\lambda_{r}(t)$:
\begin{align}
	\label{obj1}
	\Lambda_r=& \sum_{i=1}^{n_r} \int_{t_{i-1}}^{t_{i}}  \lambda_r(\xi) \cdot  \tau_r(\xi) \cdot  f_{o_r r}(\xi) \, {\rm d \xi} =  
	\sum_{i=1}^{n_r} \sum_{j=1}^{k} \left [\int_{[t_{i-1}, t_{i}] \cap [\delta_{j-1},\delta_j ]} \lambda_r(\xi) \cdot  \tau_r(\xi) \cdot  f_{o_r r}(\xi) \, {\rm d \xi}  \right ].
\end{align}

Introducing the notation 
$I_{ij}=[t_{i-1}, t_{i}] \cap [\delta_{j-1},\delta_j]$ and defining 
\[
\widehat{\lambda}_{rij}= \left \{
\begin{array}{ll}
	\lambda_{rj}& \text{ if } I_{ij} \neq \emptyset\\
	0& \text{otherwise}
\end{array}
\right .,
\]
expression (\ref{obj1}) can be rewritten as:
\begin{align}
	\label{obj2}
	\Lambda_r=&   
	\sum_{i=1}^{n_r} \sum_{j=1}^{k} \widehat{\lambda}_{rij}\int_{I_{ij}}   \tau_r(\xi) \cdot  f_{o_r r}(\xi) \, {\rm d \xi}. 
\end{align}

Finally, approximating the integral yields:
\begin{align*}
	&\Lambda_r \approx \sum_{i=1}^{n_r} \sum_{j=1}^{k} \widehat{\lambda}_{rij} \cdot  \tau_r( t_{i-1}) \cdot  f_{o_r r}(t_{i-1})  \cdot  \mu_{ij},
\end{align*}
where $\mu_{ij}$ represents the measure of the interval $I_{ij}$. Consequently, if $I_{ij}=\emptyset$, then $\mu_{ij}=0$.

On the other hand, the simulation model enables the calculation of the tons of freight transported by each mode, allowing the value of externalities to be assessed (i.e., the second term of the objective function). Adding the value associated with these externalities to the (approximated) total value of \gls{TAC} yields the objective function's value for the decision variables $\{\lambda_{rj}\}$.

\subsection{Discrete-event simulation algorithm}

In this section, the set of model constraints previously defined are discretized, specifically: system dynamics equations~(\ref{eq3})–(\ref{eq6}), the flow propagation constraints~(\ref{eq7})–(\ref{eq10}), and the modal split constraints~(\ref{demand:eq1})–(\ref{demand:eq2}). This discretization results in a discrete-event simulation model that describes the propagation of flows within the network.

The simulation model is structured as a discrete-event system, where events must be processed strictly in chronological order. To enforce this ordering, the event processing time is defined as:
\begin{equation}
    T_i \equiv \text{ the time at which event } i \text{ must be processed.}
\end{equation}
\noindent The simulation clock is given by $T = \min_i {T_i}$ and advances sequentially as each event is executed. The simulation terminates when the clock exceeds the study period, i.e., when $T \geq T_{\max}$.

Each event corresponds to a packet and represents the moment at which a packet departs from an arc. A single packet may generate multiple events, denoted as $i_1, i_2, \dots$, but at any given simulation time $T$, only one of these events is active. Thus, each packet is uniquely identified by its active event, with index $i$ referring to the corresponding active event.

Packets transition between two states: $\texttt{STATE0}$ and $\texttt{STATE1}$. In $\texttt{STATE0}$, a packet is waiting to enter the network at a certain connector, while in $\texttt{STATE1}$, the packet is actively traversing the network. The simulation mechanism differs depending on the current state of the packet.

\begin{figure}[h]
    \centerline{\includegraphics[scale=.45]{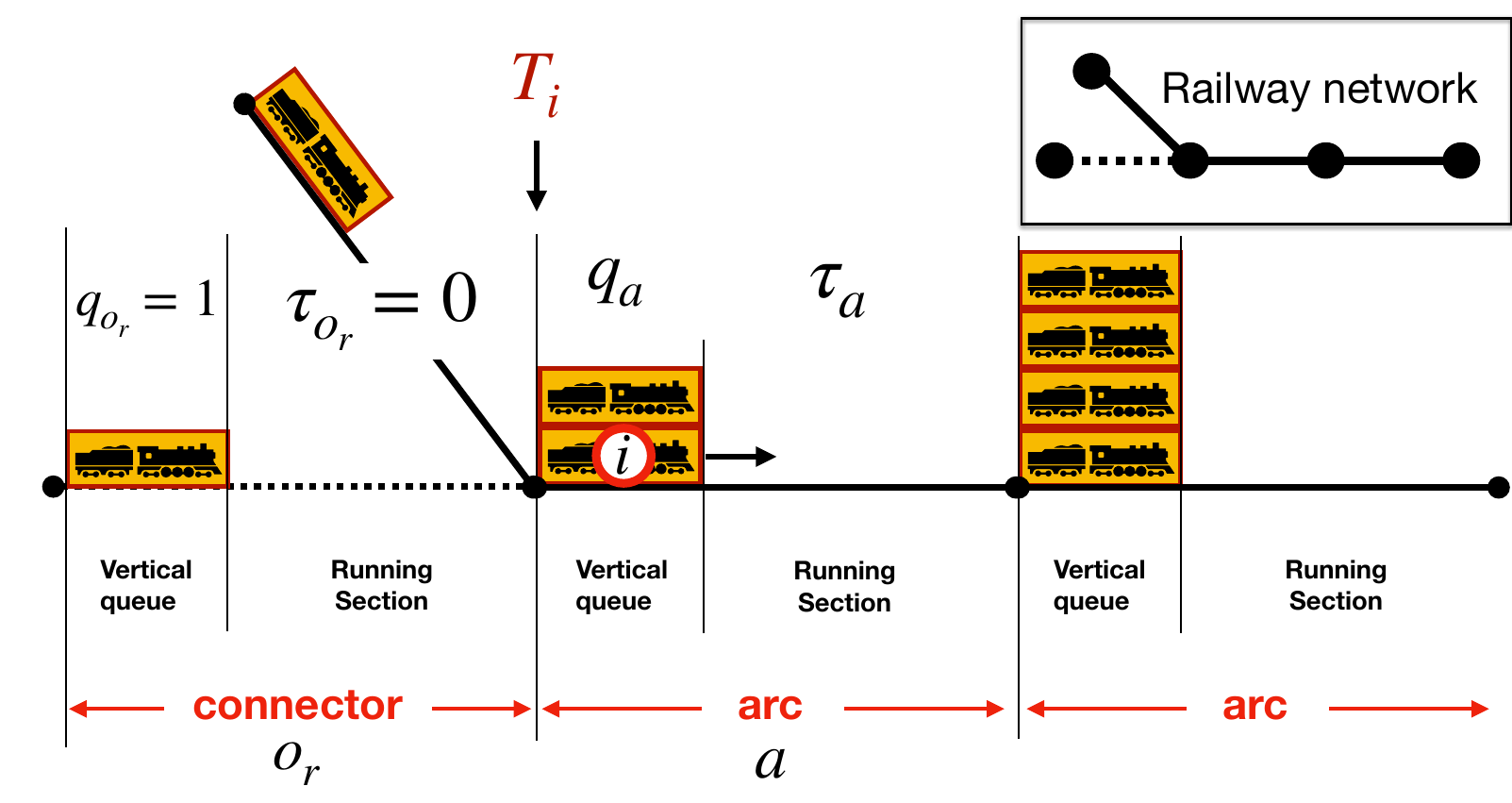}}
\caption{Representation of the railway network in the simulation model \label{fig:simulation}}
\end{figure}

Figure~\ref{fig:simulation} illustrates this process, where flow is discretized into packets of equal size $\Delta f$ (where $\Delta f = 1$ corresponds to a single train). At the current simulation time, the clock determines which packet $i$ will be processed next. A packet located on a physical arc is in $\texttt{STATE1}$ and thus enters the running section, joining the queue for its next arc. When packet $i$ enters the running section at time $T_i$, the queue becomes active, and the flow reaches its maximum capacity over the time interval $\Delta t$ required to process a packet of volume $\Delta f$:
\begin{equation}
	\phi_{ar}(t)=k_a(t), \quad \text{for all } t \in [T_i,T_i+\Delta t].
\end{equation}
Then, integrating both sides of the equation,
\begin{equation}
	\int_{T_i}^{T_i +\Delta t }\phi_{ar}(t) {\rm d}t=\int_{T_i}^{T_i +\Delta t } k_a(t) \, {\rm d}t.
\end{equation}
Under the assumption of constant capacity $k_a$ over this interval, the following is obtained:
\begin{equation}
\label{eq:tiempo_procesamiento}
	\Delta f=\int_{T_i}^{T_i +\Delta t }\phi_{ar}(t) \,  {\rm d}t=\int_{T_i}^{T_i +\Delta t } k_a \, {\rm d}t=  k_a \cdot \Delta t \Rightarrow \Delta t= \dfrac{\Delta f}{k_a}.
\end{equation}
This constraint ensures that no additional packets in the queue can enter the running section until the current packet has been fully processed. 

Therefore, the processing time for a packet is determined by Equation~(\ref{eq:tiempo_procesamiento}). Now, let $a_i$ be the current arc of packet $i$, and let $q_{a_i}$ represent the set of packets in the vertical queue of arc $a_i$.
While packet $i$ is being processed, all other packets in the queue $q_{a_i}$ must wait, which means that their events must be rescheduled to occur no earlier than the completion of the processing of packet $i$. The update rule for the next event is given by:
\begin{equation}
\label{eq:dwell_time}
 T_j= \max \left (T_j, T_i+\dfrac{\Delta f }{k_{a_i}} \right ), \, \forall j \in q_{a_i}.
\end{equation}

When packet $i$ has been processed, it updates its event time to $T_i = T_i + \tau_{a_i}$, advances to its new arc $a_i = (a_i)^+$, and enters its corresponding vertical queue. The algorithm then proceeds by processing the next event in the system.

Connectors are modeled as fictitious arcs with infinite capacity, \( k_a = \infty \), and zero travel time in the running section, \( \tau_a = 0 \). These arcs serve as entry points for generating packets (trains) that will propagate through the network. When an event \( i \) is in \( \texttt{STATE0} \), the associated packet has not yet been physically generated and remains latent. As described in Algorithm~\ref{alg:discrete_event_algorithm}, event \( i \) is processed at the moment sufficient demand has accumulated to fill the train. This instant is determined by solving an integral equation ensuring that the required demand has been reached.

At the same time, a new latent event \( j \) is generated to account for subsequent packet creation along the same path $r_i$. The timestamp of this new event, \( T_j \), coincides precisely with the moment packet \( i \) enters the network, i.e., \( T_j = T_i \).

The procedure described above is formalized in the pseudocode of Algorithm~\ref{alg:discrete_event_algorithm}, which governs the discrete event simulation. Furthermore, to ensure transparency and reproducibility, the code and experiments developed in this study are publicly available in the GitHub repository associated with this paper\footnote{\url{https://github.com/RicardoGarciaRodenas/TAC-access-charge-for-freight-rail-transport}}.

\begin{algorithm}[!h]
\DontPrintSemicolon
  \setstretch{1.1}
  \KwInput{Parametrization $\lambda_{r1}\,\ldots,\lambda_{rk}$ for all $r$.}
  \KwOutput{Network events trace.}
  \tcc{Initialization }
    $D_\omega=\int_0^{T_{\max}} D_\omega(\xi) \, {\rm d \xi}$; \tcc*{Total demand for each origin-destination pair}
    $L_\omega=D_\omega$; \tcc*{Assume that total demand is initially transported by road}
    $n_p=0$; $n_e=0$;  \tcc*{Packet and event counters}
    $E=\emptyset$; \tcc*{Event set}
    $T=0$; \tcc*{Simulation clock}
    \BlankLine
    \ForEach{$r\in \cup_\omega R_\omega$}   
    {
        $n_e=n_e+1$; $j=n_e$; $E=E\cup \{j\}$; $T_j=0$; \tcc*{Initialize a new event j}
        $state_{j}=\texttt{STATE0}$; $r_j=r$, $\omega_j=\omega$; $a_j=o_r$; $q_{a_j}=\{j\}$;\\
        
        $\tau_r=\overline{\tau}_r$; \tcc*{Initial travel time for the path $r$ equal to the commercial travel time}
        \BlankLine
    }
    \tcc{Simulation loop}
    \While{$T \le T_{\max}$}{
        $i = \underset{j \in E}{\hbox{arg minimize}} \{T_j\}$;\\
        Update the clock $T=T_i$; \\
        \BlankLine
        \If{$state_i == \texttt{STATE0}$}
        {   
           
          $n_e=n_e+1$; $j=n_e$; $E=E\cup \{j\}$;\tcc*{Generate a new event j}
          $state_{j}=\texttt{STATE0}$; $r_j=r_i$, $\omega_j=\omega_i$; $a_j=a_i$; $q_{a_j}=q_{a_j}\cup\{j\}$;\\
        \tcc{Compute the instant associated with the event $j$}
        \ForAll{$r\in R_{\omega_i}$}
            {$U_{r}= \left ( A_{r}(T_i) +  \lambda_r(T_i) \right ) \cdot \tau_{r} +C_{r}$; \tcc*{Parameterized by the \glspl{TAC} $\lambda_{r1},\ldots \lambda_{rj}$}  
            }
        Compute the instant $T_j$ as the solution to  the integral equation: \\
        $\;\;\;\;\;\;\;\;\Delta f  = \kappa\cdot \left [\dfrac{\exp \left (U_{r_i} \right ) }{\exp(V_{\omega_i})+\sum_{r' \in R_{\omega_i}}\exp \left (U_{r'}\right )} \right ] \cdot
        \int_{T_{i}}^{T_j}  \cdot D_{\omega_i}(t) \, {\rm dt};$
        \BlankLine
        $state_i=\texttt{STATE1}$; $n_p=n_p+1$; $T_i=T_j$;  $t^0_i=T_i$; \tcc{Load packet $i$ onto the network} 
        $L_{\omega_i}=L_{\omega_i}-\frac{\Delta f}{\kappa}$; \tcc{Update the freight for road}
            }
        \If{$a_i == d_{r_i}$}{
            $E=E - \{i\}$; $\tau_{r_i}=T_i-t^0_i$;  \tcc*{The packet  has reached its destination}
        }
        \Else{
            $T_j= \max \left (T_j, T_i+\dfrac{\Delta f }{k_{a_i}} \right ), \forall j \in q_{a_i};$ \tcc*{Update entry time to the running section}
            $T_i=T_i+\tau_{a_i}$; $q_{a_i}=q_{a_i} -\{i\}$; \tcc*{Packet i crosses  the running section}
            $a_i=(a_i)^+$; \tcc*{Obtain the next arc}
            $q_{a_i}=q_{a_i}\cup \{i\}$; \tcc*{Queue the packet}
        }
    }
\caption{Discrete Event Simulation Algorithm}
 \label{alg:discrete_event_algorithm}
\end{algorithm}

\section{Numerical experiments}

This section presents a series of numerical experiments carried out using the proposed model. Specifically, the Mediterranean Rail Freight Corridor is analyzed due to its particular importance in Europe. The study first evaluates the feasibility of solving the proposed model, followed by an analysis of the model's outputs for the case study.

\subsection{The Mediterranean Rail Freight Corridor (RFC6)}

\begin{figure}[h]
    \centerline{\includegraphics[scale=.7]{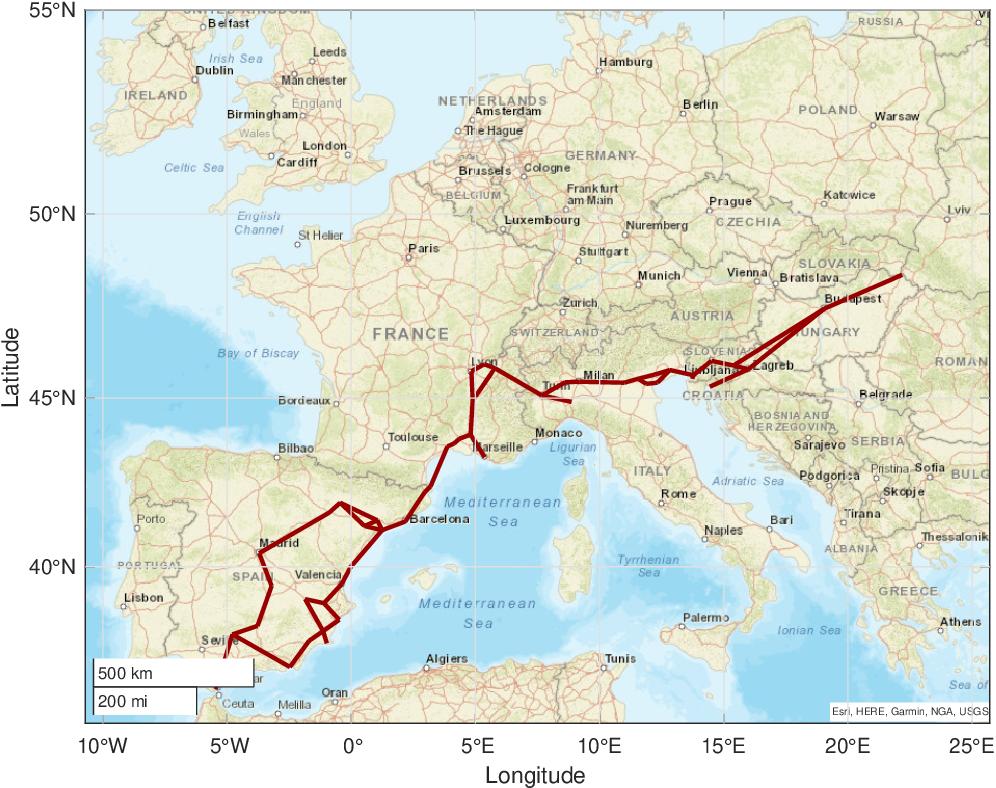}}
\caption{The Mediterranean Rail Freight Corridor (RFC6)\label{fig:CM}}
\end{figure}

The Mediterranean Corridor (see Figure~\ref{fig:CM}) stands out as Europe's primary horizontal rail freight axis, extending $7967$ kilometers from Spain to the European Union border. It serves as a critical link between the Mediterranean Basin and Central Europe, including Ukraine. Running East to West, the corridor intersects seven other Freight Railway Corridors and traverses three of Europe's four major manufacturing regions: Catalonia, Auvergne-Rhône-Alpes, and Piedmont-Lombardy. The countries along the corridor collectively account for a GDP of approximately $5800$ billion euros and a population of around $190$ million. Additionally, the corridor facilitates connectivity with over $100$ intermodal terminals and provides direct links to five major Mediterranean seaports as well as two key river ports in Lyon and Budapest. 

Figure~\ref{fig:CM} illustrates the routes of the corridor: Almería-Valencia / Algeciras / Madrid-Zaragoza / Barcelona-Marseille-Lyon-Turin-Milan-Verona-Padua / Venice-Trieste / Capodistria-Ljubljana-Budapest-Ljubljana / Fiume-Zagabria-Budapest-Zahony (Hungarian-Ukranian border).

This study assumes a centralized authority responsible for setting \glspl{TAC}, integrating all relevant infrastructure managers, including: ADIF (Spain), Línea Figueras Perpignan (Spain-France), SNCF Réseau (France), Oc'Via (France), RFI (Italy), SŽ - Infrastruktura (Slovenia), MÁV (Hungary); VPE (Hungary), and HŽ Infrastruktura (Croatia).

\subsubsection{Origin-destination freight flows}

This section presents the estimation of freight flows along the Mediterranean Corridor. The base information is the corridor information documents \citep{RFC6-21a,RFC6-21b} titled {\sl Mediterranean RFC Implementation Plan TT (2021/2022)\footnote{\url{https://www.medrfc.eu/wp -content/uploads/2022/01/rfc6_implementation-plan-tt-2022-complete-07-01-2022.pdf}} and (2022/2023)\footnote{\url{https://www.medrfc.eu/ wp-content/uploads/2022/01/2-med-rfc-ip-tt-2023-complete-11-01-2022.pdf}}}. These documents include a market analysis study and provide projections for $2030$, considering various GDP and transport cost evolution scenarios. The projections assume that rail competitiveness will align with expectations if the corridor is fully implemented by $2030$. The scenarios considered are named {\sl Scenario $2$} and {\sl Scenario $3$}, which evaluate the trend macroeconomic case combined with different transport cost evolution assumptions (worst and best cases).

The scenarios define rail traffic volumes and modal shares for country-to-country relations. In addition to internal freight flows, the study estimates exchange and transit flows, incorporating not only the corridor countries but also broader aggregated regions, including South-Eastern Europe, North-Eastern Europe, and Western Europe. To adapt this data for the proposed model, a disaggregation of flows at the node level is necessary. Internal flows are distributed proportionally to the population of each node, while transit and exchange flows are assigned to ports, nodes with terminals, and nodes connected to other railway corridors, projecting these flows proportionally to the number of terminals per node or the volume handled by the port.

The railway network under consideration consists of $49$ nodes, some of which correspond to key junctions essential for accurately defining railway segments. Finally, the estimated freight flows for the $612$ main O-D pairs are obtained. The flow intensity is assumed to remain constant throughout the planning period, representing a full year, i.e., \( D_\omega(t) = d_\omega \). This constant is calibrated to ensure the model accurately replicates the transport demand forecast for {\sl Scenario 3} in the year 2030.

\begin{figure}[h]
    \centerline{\includegraphics[scale=.7]{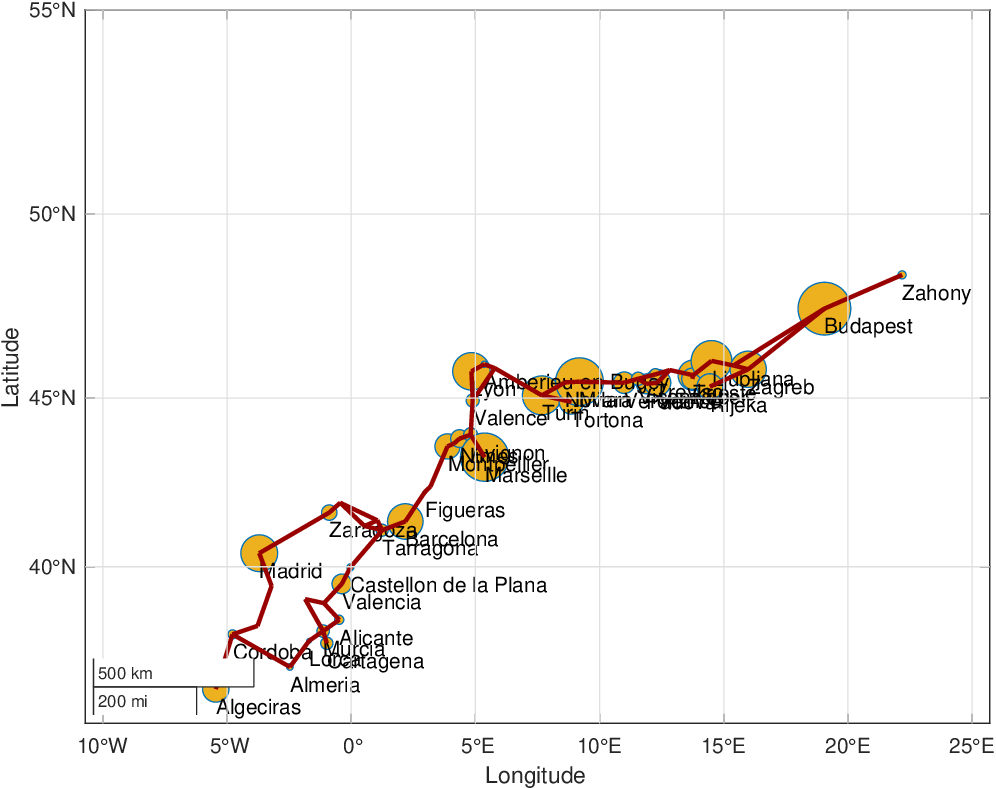}}
\caption{Desegregation of freight demand by nodes\label{fig:demand_nodes}}
\end{figure}

Figure~\ref{fig:demand_nodes} illustrates the significance of each node within the estimated O-D flow matrix, highlighting their respective roles as major generators or attractors of freight along the corridor.

\subsubsection{The logit model: market shares}

The analysis assumes that for each origin-destination pair, only the minimum-cost path is used, implying that $|R_\omega|=1$. A static logit model serves as the foundation for developing a dynamic one. For a given origin-destination pair $\omega = (o, d)$, the choice between road and rail transportation is modeled using two distinct utility functions. The rail market share is determined by:
\begin{equation}
    Logit(V_\omega,U_\omega)=\dfrac{\exp(U_\omega)}{\exp(V_\omega)+\exp(U_\omega)},
\end{equation}
\noindent where $U_\omega$ and $V_\omega$ represent the utilities associated with rail and road transport, respectively. The utility functions are assumed to follow a linear form:
\begin{align}
V_\omega&= \beta_{road} C_{road} + \alpha_o+\alpha_d,\\
U_\omega&= \beta_{railway} C_{railway},
\end{align}
\noindent where the parameters $\alpha_o$ and $\alpha_d$ are associated with the origin and destination nodes, respectively. To simplify the estimation, all nodes within the same country share a common parameter $\alpha$. If a node functions as both an origin ($o$) and a destination ($d$), its parameter satisfies $\alpha_o = \alpha_d$. To maintain consistency across all origin-destination pairs, costs are expressed in terms of cost per kilometer per tonne \(\left( \text{\euro} \, / \text{ton} \cdot \text{km} \right)\).

The total cost of railway transportation is defined as:
\begin{equation}
C_{railway}=  c_T+c_{AccessCharges}+c_\ell,
\end{equation}
where  $c_T$ represents the delay cost \(\left( \text{\euro} \, / \text{ton} \cdot \text{km} \right)\),  $c_{AccessCharges}$ accounts for the network access charges \(\left( \text{\euro} \, / \text{ton} \cdot \text{km} \right)\), and $c_\ell$ denotes the fixed cost \(\left( \text{\euro} \, / \text{ton} \cdot \text{km} \right)\) when traveling at the standard commercial speed $\overline{v}_r$ on the path $r$. Note that the fixed component of the train path cost is $C_r = c_\ell \ell_r$, where $\ell_r$ is the length (in kilometers) of path $r$.

The delay cost quantifies the total additional expense incurred by the railway operator when a train travels at a speed lower than the reference commercial speed for the route, denoted as \( \bar{v}_r \). This cost is computed as \( c_t( \tau_r(t) - \overline{\tau}_r ) \), where \( c_t \) is the cost per tonne and hour of delay \(\left( \text{\euro} \, / \text{ton} \cdot \text{h} \right)\), and \( \bar{\tau}_r \) represents the travel time required to traverse route \( r \) at the reference commercial speed \( \bar{v}_r \). Since \( c_T \) represents the cost per kilometer and tonne, the cost per unit distance is obtained by dividing by the length of route \( r \), \( \ell_r \), i.e.,
\begin{equation}
c_T = \dfrac{
c_t( \tau_r(t) -  \overline{\tau}_r )}
{\ell_r} =  \frac{c_t}{\ell_r} \left(1 - \frac{\overline{\tau}_r}{\tau_r(t)}\right) \cdot \tau_r(t)=
\underset{A_r(t)}{\underbrace{\frac{c_t}{\ell_r} \left( 1-  \frac{v_r (t)}{\overline{v}_r} \right)} }\cdot \tau_r(t),
\end{equation}
where \( v_r(t) \) denotes the average travel speed on path \( r \) when the journey begins at time \( t \), given by  
$
v_r(t) = \frac{\ell_r}{\tau_r(t)}. 
$ In this equation, the term \( A_r(t) \) identifies this cost within the general cost structure expression.

The \glspl{TAC} are modeled as a proportion of the railway operator's fixed costs:
\begin{align}
    \lambda_r(t)&= \frac{p_r(t) }{\overline{\tau}_r}\cdot c_\ell,\\
    c_{AccessCharges}&= \lambda_r(t) \cdot \tau_r(t),
\end{align}
\noindent where it is assumed that $0\le p_r(t) \le 0.25$, that is, that the fees cannot represent more than $25\%$ of the reference costs of the rail operator. This upper bound is established a priori as a regulatory constraint.

The cost for both rail transport ($c_\ell$) and road transport ($C_{road}$) are based on the values reported in \cite{SWW23}. As observed in \cite{RFC6-21b}, commercial train speeds vary across different segments of the corridor. For this study, an average commercial speed is adopted for the entire network to maintain consistency. The parameters used in the cost model are summarized in Table~\ref{tab:cost_parameters}.

\begin{table}[h]
\centering
\caption{Parameters associated with the transportation costs}
\label{tab:cost_parameters}
\resizebox{0.8\textwidth}{!}{
\begin{tabular}{cccc}
\toprule
 $c_t$ (\euro \, / ton $\cdot$ h) & $c_\ell$ (\euro  \, / ton $\cdot$ km) & $C_{road}$ (\euro \, / ton $\cdot$  km) & $\overline{v} (km/h)$\\
\hline
2.2300& 0.0450 & 0.3850 & 53\\
\bottomrule
\end{tabular}
}
\end{table}

To calibrate the logit model parameters, scenarios $2$ and $3$ from \cite{RFC6-21a} are used. The primary difference between these scenarios is an $18\%$ increase in road transportation costs. Using the specified costs and considering each entry in the railway market share matrix for each scenario as an observation of the model, the logit model is adjusted accordingly. The estimated parameters are presented in Table~\ref{tab:logit_parameter}.

\begin{table}[h]
\centering
\caption{Parameters associated with the logit model \label{tab:logit_parameter}}
\label{tab:LogitModel}
\resizebox{0.8\textwidth}{!}{
\begin{tabular}{cccccccc}
\toprule
$\alpha_{ES}$& $\alpha_{FR}$ &  $\alpha_{IT}$ &  $\alpha_{SL}$  & $\alpha_{HK}$  & $\alpha_{HU}$  & $\beta_{railway}$ & $\beta_{road}$ \\
\hline
0.5520 & 0.4589 & 0.1356 & 0.3512 &0.2220 &0 & -149.8372 & -13.5454 \\
\bottomrule
\end{tabular}
}
\\ES (Spain) FR (France) IT (Italy) SL (Slovenia)  HK (Croatia) HU (Hungary)
\end{table}

\subsubsection{Policies for track access charges}

The pricing of railway network access is established based on both economic efficiency and environmental-road safety considerations. These criteria are incorporated into the objective function through the weighting parameter $\eta$. In this section, three pricing policies are introduced, each defined by specific values of $\eta$, which will be analyzed below.

The emission rate for a default diesel train used in these experiments is set at $23$ gCO$_2/$tonne-km. This value is consistent with the parameters obtained from the five models tested in \cite{Hei20} and aligns with the reference parameter used in the Network for Transport Measures model\footnote{\url{https://www.ecotransit.org/en/emissioncalculator/}}. For road transport emissions, significant variability exists depending on vehicle characteristics. For instance, one study\footnote{\url{https://theicct.org/publication/co2-emissions-from-trucks-in-the-eu-an-analysis-of-the-heavy-duty-co2-standards-baseline-data/}} found that urban delivery trucks with a 4x2 axle configuration (4-UD) emit an average of $307$ gCO$_2/$tonne-km, which is more than five times higher than the emissions of long-haul tractor-trailers (5-LH), which average $57$ gCO$_2/$tonne-km.

To estimate emissions for both road and rail transport, the EcoTransIT World simulator\footnote{\url{https://www.ecotransit.org/en/emissioncalculator/}} was used. The longest segment of the Mediterranean Corridor, Algeciras to Záhony (3277 km), was considered. The simulation verified that truck emissions are between $2.4$ and $6.5$ times higher than those of rail transport (whether using electric or diesel locomotives). These findings align with projections by \cite{JVH24} for freight transport emissions in Sweden, which estimate that by 2040, assuming a fuel mix with $70\%$ biofuels, road transport emissions per tonne-kilometer would still be 5.5 times higher than those of rail transport.

Based on these findings, the ratios $2.4$ and $6.5$ are used as reference values in this study. Using the train emission rate of $23$ gCO$_2/$ton-km, the corresponding truck emission rates are set at $54.0$ gCO$_2$/ton-km and $149.7$ gCO$_2$/ton-km, respectively.


The environmental benefit of rail transport is quantified by converting CO$_2$ emissions reductions into monetary values using carbon credits. A carbon credit represents one tonne of CO$_2$ equivalent (CO$_2$e) that an entity is permitted to emit. In this study, the carbon price is set at $54.21$\euro{} per tonne of CO$_2$e, based on market values as of February 25, 2024.

The parameter $\eta$ represents the monetary value of CO$_2$e reductions achieved by shifting freight transport to rail. Three policies are defined based on different emission reduction assumptions:
\begin{align*}
    \hbox{\bf Policy 1: } &   \eta_{1}=\dfrac{149.7-23}{10^6} \times \text{creditValue}\\
    \hbox{\bf Policy 2: } &   \eta_{2}=\dfrac{54.0-23}{10^6} \times \text{creditValue}\\
    \hbox{\bf Policy 3: } &   \eta_3=0
\end{align*}

\subsection{Resolution of the model}
\label{sec:exp_comp}

Several challenges arise when trying to solve the proposed optimization model. Firstly, evaluating the objective function is computationally expensive, as it relies on a simulation model that traces the movement of approximately $32,500$ trains per year, transporting around $40$ million tons of freight. Secondly, the objective function is non-differentiable due to the discrete-event nature of the simulation. In this context, an infinitesimal increase in the train-path fee does not necessarily lead to changes in the simulation output, making the application of exact gradient-based optimization algorithms particularly challenging.

To address these difficulties, both exact and metaheuristic optimization algorithms were tested. For exact algorithms, MATLAB implementations of the \gls{AS}, \gls{IP}, and \gls{SQP} methods were used\footnote{The {\tt fmincon} function was used with the algorithms {\tt active-set}, {\tt interior-point}, and {\tt sqp}.}. For metaheuristic approaches, MATLAB implementations of \gls{PSO}, \gls{GA}, and \gls{PS} were applied\footnote{The functions {\tt particleswarm}, {\tt ga}, and {\tt patternsearch} were used.}.

Two distinct pricing schemes were assessed:
\begin{enumerate}
    \item {\bf Proportional approach.} In this pricing scheme, access charges are applied proportionally to the length of each railway path. A single proportionality constant is used for all routes, ensuring uniform pricing across the network. This constant is defined as a fraction of the fixed rail transportation cost, leading to a toll structure expressed as:
    \[
    p_r(t) = p, \quad \text{ for all }r, \text{ for all } t.
    \]  
    This approach is equivalent to optimally determining a conventional per-kilometer toll scheme.
        
	\item {\bf Path-based approach.} In this scheme, each railway path is assigned a distinct proportionality constant, allowing for differential pricing based on specific characteristics of individual routes:
    \[
    p_r(t) = p_r, \quad  \text{ for all } r.
    \]  
\end{enumerate}

The proportional approach offers a straightforward practical implementation and a simple solution, as it reduces the optimization problem to a single variable ($p$). In contrast, the path-based approach is more general and is expected to yield better results, but it also introduces greater complexity into both its resolution and practical implementation.

\begin{figure}[!h]
	\begin{center}
		\includegraphics[width=0.32\textwidth, angle=0, origin=c]{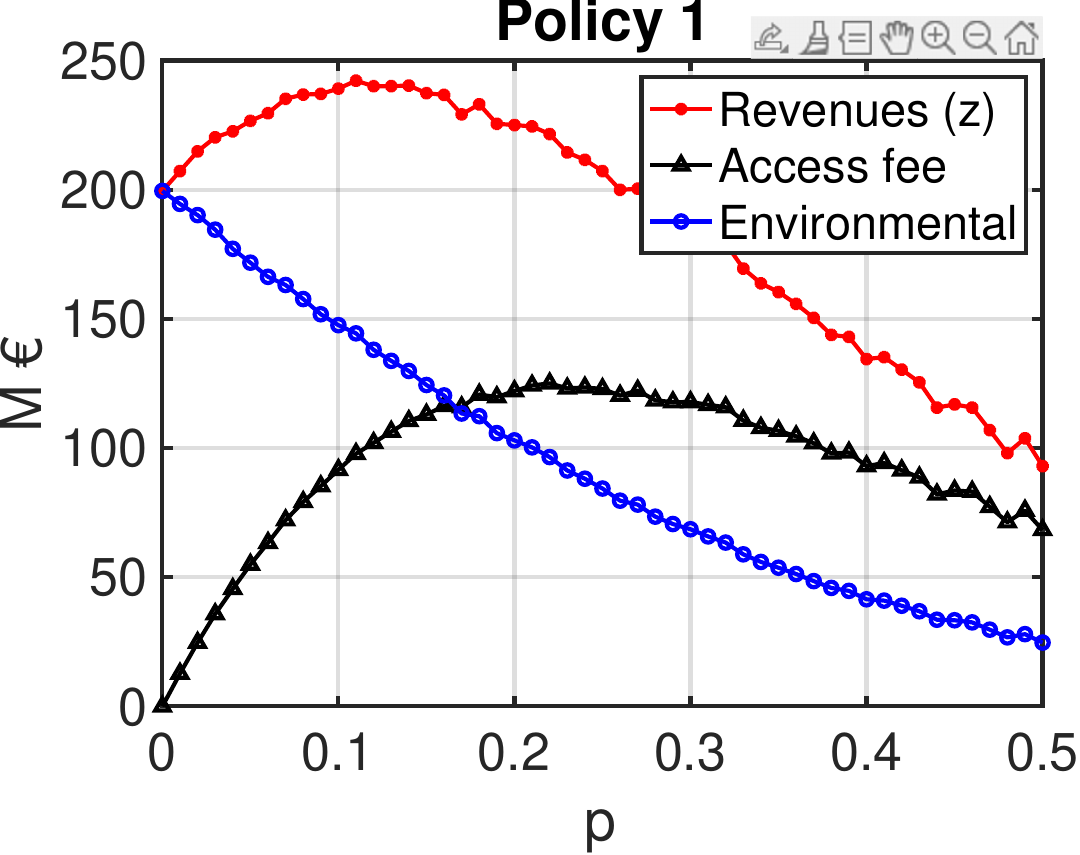}
		\includegraphics[width=0.32\textwidth, angle=0, origin=c]{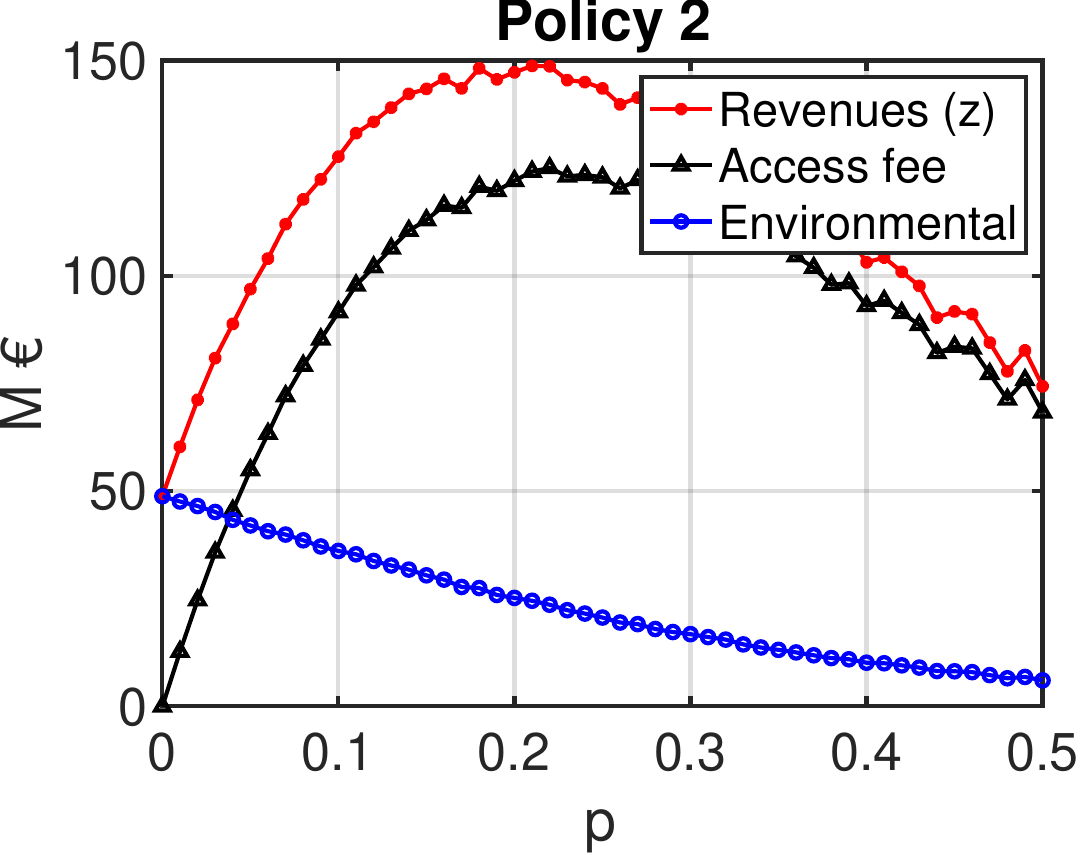}
        \includegraphics[width=0.32\textwidth, angle=0, origin=c]{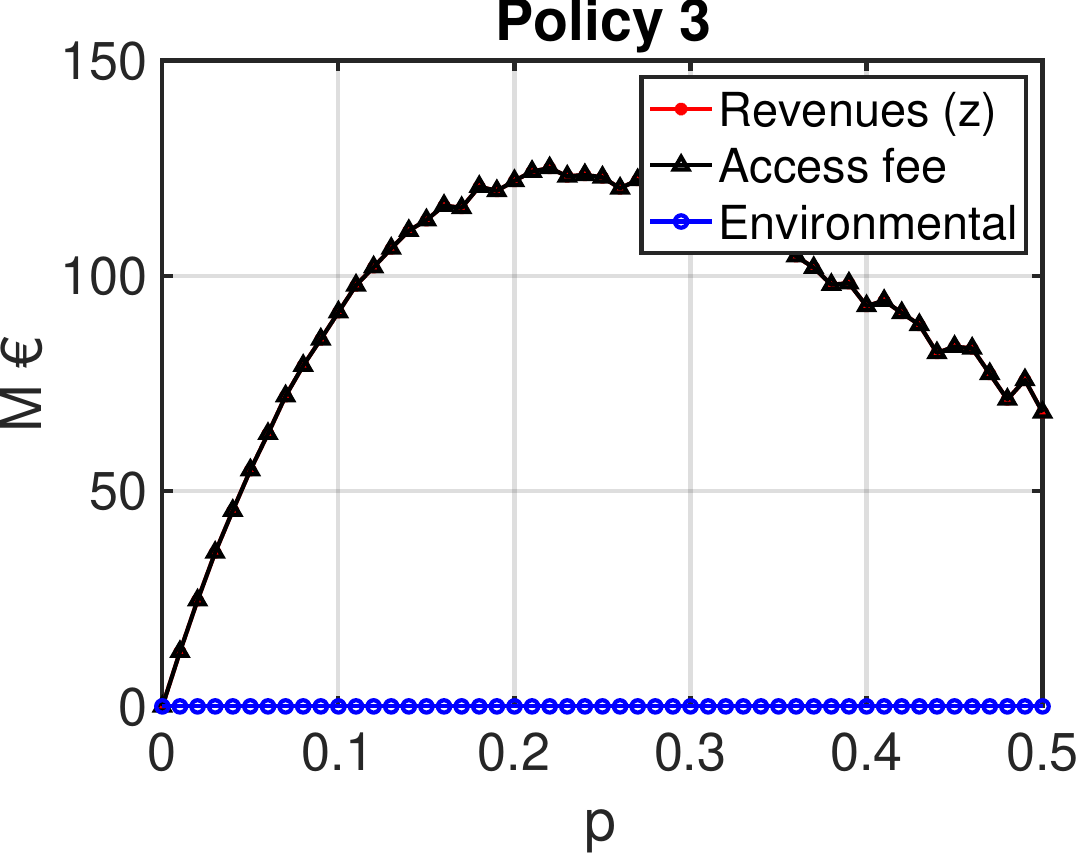}

	\end{center}
 \caption{Objective function for the proportional approach\label{fig:proportional_approach}}
\end{figure}

Figure~\ref{fig:proportional_approach} illustrates the objective function under the proportional pricing scheme. The objective function, shown in red, represents the sum of Access Charges (in black) and Environmental Costs (in blue). In Policy 3, where environmental costs are not considered, the objective function is represented solely in black. This model can be solved using a brute-force approach, where the objective function is evaluated over uniformly distributed values of $p$ on the interval $[0, 0.25]$. The figure highlights the non-differentiability of the objective function with respect to $p$, a consequence of the fact that the calculation of economic indicators disregards the number of trains circulating in the network at the end of the simulation period (1 year).

To solve the path-based approach, optimization algorithms are required. The algorithms \gls{AS}, \gls{IP}, \gls{SQP}, and \gls{PS} are initialized from an initial solution. The solution obtained from the proportional approach is used as the starting point, providing a solid initial solution for these methods to refine and improve.

A strategy proven effective in the field of transportation optimization \citep{EGR15,LGG14,GLS19,ACG14} is the hybridization of metaheuristic and exact algorithms. A simple hybridization approach applied to \gls{GA} and \gls{PSO} involves refining the solution obtained by the metaheuristic algorithm using the \gls{AS} method, starting from the output of the metaheuristic as the initial solution. This improvement has been incorporated into the metaheuristic algorithms employed in this study.

The computations were performed using MATLAB R2023b on a Mac computer, equipped with an Apple M3 Pro processor, 18GB of RAM, and running macOS Sonoma 14.3.1. For the population-based algorithms, the option {\tt PopulationSize} was set to 60, and {\tt MaxGenerations} to 500. For the exact algorithms, the maximum number of function evaluations was set to $4000$ ({\tt MaxFunctionEvaluations}). All algorithms met these termination criteria without being stopped prematurely, except for \gls{PS}, for which the number of function evaluations was extended up to $75,000$. Additionally, for the algorithms \gls{GA}, \gls{PSO}, and \gls{PS}, parallel evaluations of the objective function were implemented using $12$ workers to enhance computational efficiency.

The results are presented in Table~\ref{tab:ComputationalCost}. The first notable observation is that, in many cases, the exact algorithms do not improve upon the solution provided by the proportional approach. For instance, in {\sl Policy 1}, the \gls{IP} and \gls{SQP} algorithms fail to outperform the proportional solution, which has a value of $241.41$ M\euro. Additionally, population-based metaheuristic algorithms also do not surpass the optimal solution achieved by the proportional approach.

These findings suggest that the most effective algorithm for this problem is \gls{PS}. However, its computational cost is considerably higher than that of other methods, requiring over $10$ hours of runtime. This increased computational effort is due to the fact that \gls{PS} relies entirely on function evaluations and systematically explores the feasible region, effectively addressing the insensitivity of the simulation model to infinitesimal changes in the optimization variables.

\begin{table}[h]
\centering
\caption{Computational costs}
\label{tab:ComputationalCost}
\resizebox{\textwidth}{!}{
\begin{tabular}{llllllllc}
\toprule
& \multicolumn{7}{c}{\bf Path-based apporach}& \multirow{ 2}{*}{
\begin{minipage}{2cm}
\bf Proportional approach
\end{minipage}
}
\\
\cline{3-8} 
& & \bf \gls{AS} & \bf \gls{IP} & \bf  \gls{SQP} & \bf \gls{PSO} & \bf \gls{GA} & \bf \gls{PS} &  \\ 
\hline
\multirow{ 2}{*}{\bf Policy 1}&Z (M \euro) & 242.71 & 241.41 & 241.41 & 237.88 & 242.52 & 246.24 & 241.41 \\ 
&CPU(s) & 9992.1 & 2447.5 & 2479.8 & 13173.5 & 7683.6 & 37864.4  & 54.5\\ \hline
\multirow{ 2}{*}{\bf Policy 2}&Z (M \euro) &  150.32 & 150.14 & 149.21 & 143.97 & 145.35 & 151.76 & 149.21 \\ 
&CPU(s) &13478.2 & 5671.5 & 1964.2 & 8167 & 8133.1 & 30318.8 & 53.7 \\ 
\hline
\multirow{ 2}{*}{\bf Policy 3}&Z (M \euro) & 127.11 & 119.94 & 125.27 & 116.35 & 112.8 & 128.38 & 125.27 \\ 
& CPU(s) & 8426.9 & 12641.3 & 1672.4 & 18913.4 & 4736.9 & 26135.4 & 38.7 \\ 
\bottomrule
\end{tabular}
}
\end{table}

\subsection{Experiment I: Establishing track access charges from the \gls{IM}'s perspective}

This section illustrates the application of the proposed methodology for determining \glspl{TAC}, with the year of analysis set to 2030. The first step involves model calibration, which carried out using {\sl Scenarios} $2$ and $3$ from \cite{RFC6-21a}, representing the worst-case and best-case conditions, respectively. However, the study referenced also considers a third scenario for $2030$, referred to as {\sl Scenario $1$}, which is considered the most focused projection. Under this scenario, the estimated freight volume is $232.8$ million tons, distributed as follows: $55.2$ million tons transported by rail, accounting for $23.7\%$ market share; $144.0$ million tons ($61.8\%$) transported by road; and $24.7$ million tons ($10.63\%$) moved via Short Sea Shipping.

To ensure that the base model yields a result consistent with {\sl Scenario $1$}, transportation costs are assumed to be $18\%$ higher than those reported in Table~\ref{tab:cost_parameters}. With these adjustments, the model outputs indicate a rail freight demand of $56.6$ million tons, corresponding to a $24.33\%$ market share.

The analysis considers the three  policies introduced in Section~\ref{sec:exp_comp} and calculates the fees using both the proportional method and the path-based approach. The results are summarized in Table~\ref{tab:experiment_1_costs}. The first column displays the revenue from \glspl{TAC} received by the \gls{IM}; the second shows the value of emissions reductions; the third details the cost of transporting goods by rail; the fourth indicates the cost incurred by the \gls{FOC} due to variations in commercial speed relative to the reference speed; the fifth presents the new (average) commercial speed of the rail network; the sixth shows the millions of tons transported by rail; and finally, the last column shows the corresponding rail market share.

\begin{table}[h]
\centering
\caption{Summary of numerical results for different pricing schemes\label{tab:experiment_1_costs}}
\label{tab:Results}
\resizebox{\textwidth}{!}{
\begin{tabular}{llp{1.2cm}p{1.2cm}p{1.2cm}p{1.2cm}p{1.2cm}p{1.2cm}p{1.2cm}}
\toprule
&& \bf \glspl{TAC} (M \euro) & \bf CO2e rights  (M \euro )& \bf  Trans. cost (M \euro ) & \bf  Delay cost (M \euro)& \bf Speed Aver. (km/h)& \bf Tons (M) & \bf Rail Share ($\%$) \\ 
\toprule
\multirow{ 2}{*}{\bf Policy 1}& Path-based & 105.51 & 140.74 & 922.08 & -88.33 & 62.38 & 39.74 & 17.45 \\ 
& Proportional & 101.69 & 140.26 & 918.95 & -88.9 & 62.58 & 40.08 & 17.6 \\ 
\hline 
\multirow{ 2}{*}{\bf Policy 2}& Path-based & 123.39 & 28.36 & 759.53 & -88.5 & 64.08 & 33.72 & 14.81 \\ 
&Proportional & 121.4 & 27.8 & 744.54 & -85.75 & 64 & 33.04 & 14.51 \\ 
\hline 
\multirow{ 2}{*}{\bf Policy 3}&Path-based & 128.38 & 0 & 666.98 & -84.85 & 65.1 & 30.03 & 13.19 \\ 
&Proportional & 125.27 & 0 & 645.69 & -80.8 & 65.1 & 29.02 & 12.75 \\ 
\bottomrule
\end{tabular}
}
\end{table}

The first key observation is that the path-based approach generates $2\%–4\%$ more revenue from access fees for the \gls{IM} compared to the proportional scheme, although its implementation within the railway network and market is more complex. The second important observation underscores the necessity of incorporating environmental costs into the analysis. In the environmental policies, namely {\sl Policy 2} and {\sl Policy 1}, it is evident that the reduction in revenue for the \gls{IM} due to lower access fees is fully compensated by the social benefits derived from reduced environmental costs. Specifically, the return rate is calculated as 
$
\frac{27.8}{125.27 - 121.4} = 7.18
$
in {\sl Policy 2} and  
$
\frac{140.26}{125.27 - 101.69} = 5.94
$
in {\sl Policy 1}. Given that \glspl{IM} are public entities, these social benefits should be explicitly considered in their decision-making processes.

Another significant consideration for \glspl{IM} is railway capacity management. In this numerical analysis, it is assumed that the railway system is shared with passenger-rail transport, with priority given to passenger services during daytime slots. In $2030$, the corridor is expected to have all sections as two-way tracks, with some sections featuring $n_a$ tracks per direction. The maximum capacity for a segment $a$ is defined as $k_a = 6 n_a$ vehicles per hour. Given that passenger trains primarily operate during daylight hours, the residual capacity available for freight trains is significantly reduced during these periods. To determine the available capacity at a given time instant $t$, denoted by $k_a(t)$, the time period within the 24-hour cycle must be identified using the modulo operation \( \mod(t, 24) \), which computes the remainder when time \( t \) (in hours) is divided by 24. The capacity of railway segments imposes a minimum dwell time between packets, see Equation~(\ref{eq:dwell_time}), calculated as follows:
\begin{equation}
\label{eq:capcity_residual}
T_{q_a} (t)= \left \{
\begin{array}{ll}
 \frac{\Delta f}{k_a* 0.15}  & \hbox{if } mod(t,24) \in [10,18), \\
   \frac{\Delta f}{k_a* 0.30}  & \hbox{if }   mod(t,24) \in [7,10) \cup [18,24), \\
     \frac{\Delta f}{k_a} &  \hbox{otherwise.}
\end{array}
\right .
\end{equation}

\begin{figure}[h]
	\begin{center}
		\includegraphics[width=0.32\textwidth, angle=0, origin=c]{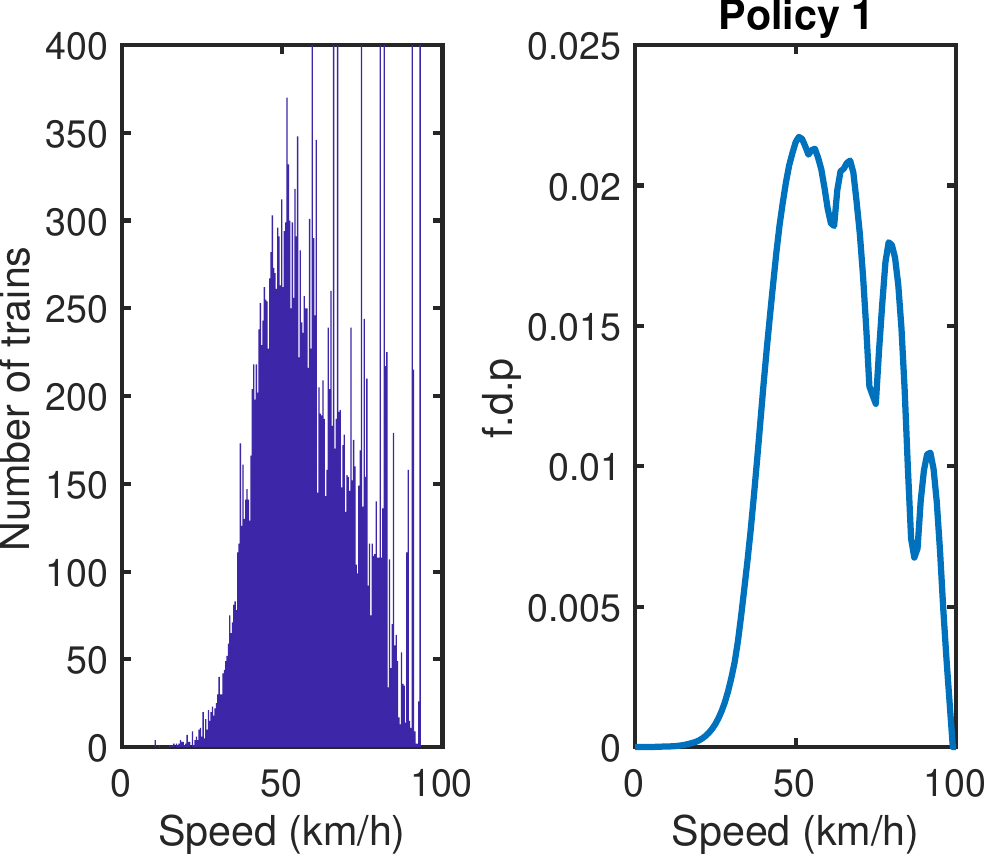}
		\includegraphics[width=0.32\textwidth, angle=0, origin=c]{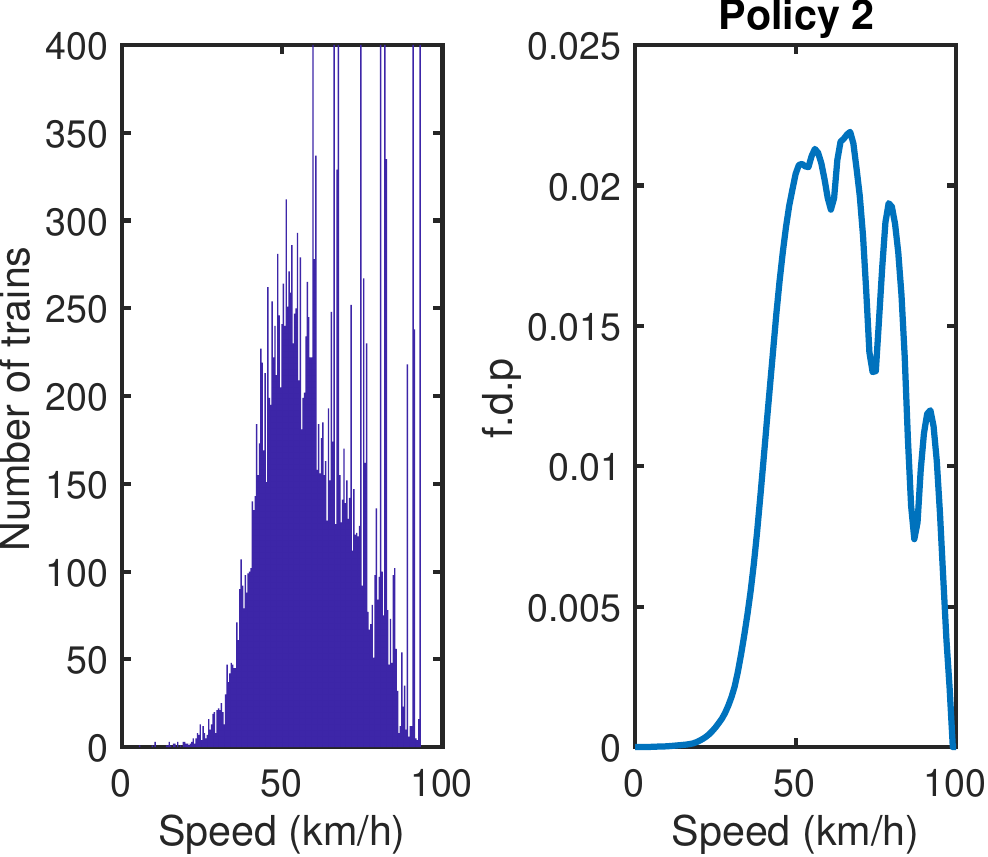}
  \includegraphics[width=0.32\textwidth, angle=0, origin=c]{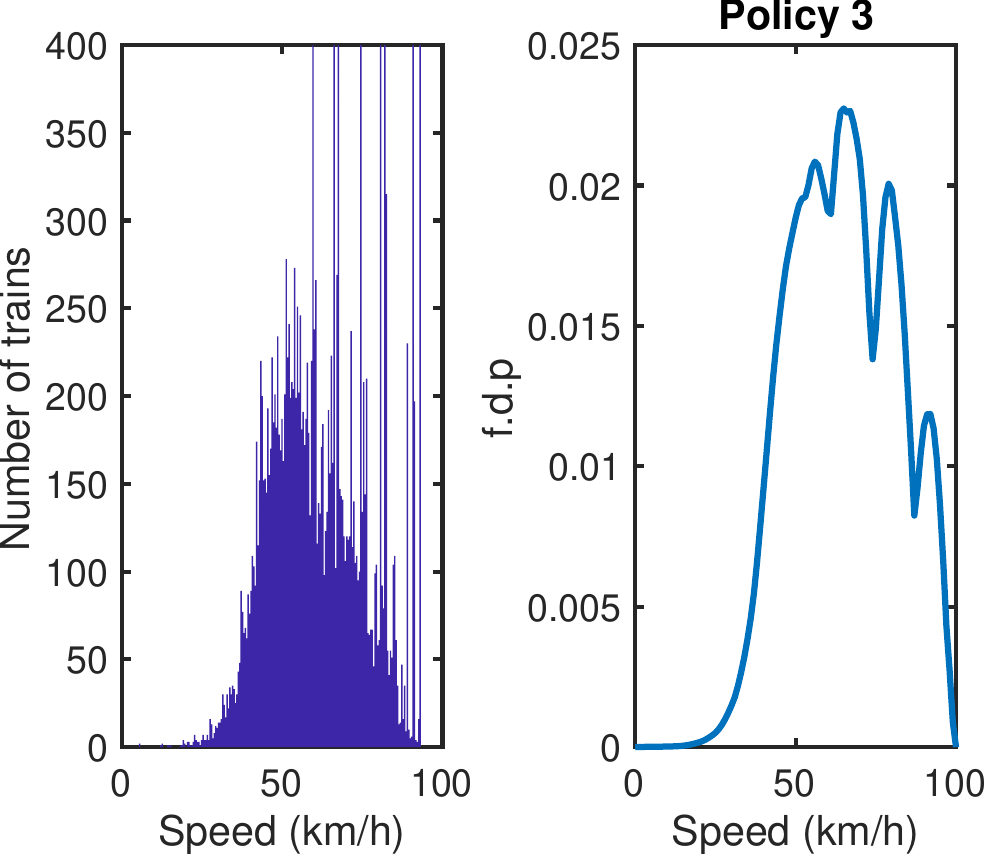}
	\end{center}
 \caption{Speed profiles across different policies \label{fig:speed_profile}}
\end{figure}

It is assumed that trains operate at a speed of $100$ km/h in the running sections. Track capacity constraints, particularly in the vertical queue segments, cause a reduction in the average speed. Currently, the commercial speed is approximately $53$ km/h, which has been used as a reference value of $\bar{v}_r$ for all paths $r$ when calculating delay costs. The model has the capability to estimate the commercial speed of trains, which varies depending on the level of network congestion. This parameter is shown in the fifth column of Table~\ref{tab:experiment_1_costs}. It is observed that the current commercial speed increases from $53$ km/h to approximately $64$ km/h, resulting in operating-cost reductions for the \gls{FOC} ranging from $80.8-88.9$ M \euro. Another important observation is that, although the increase in market share is significant, its effect on network congestion (and consequently on commercial speed variation) is minimal. This is clearly illustrated in Figure~\ref{fig:speed_profile}, which shows the speed profile of trains across the three policies. 

Another important feature of the proposed model is its ability to evaluate the elasticity of demand in response to \glspl{TAC}. Differences of up to $25\%$ in the amount of tonnage transported are observed between {\sl Policy 1} and {\sl Policy 3}, highlighting the sensitivity of transported tons to pricing structures.

\begin{figure}[!h]
	\begin{center}
 \includegraphics[width=0.9\textwidth, angle=0, origin=c]{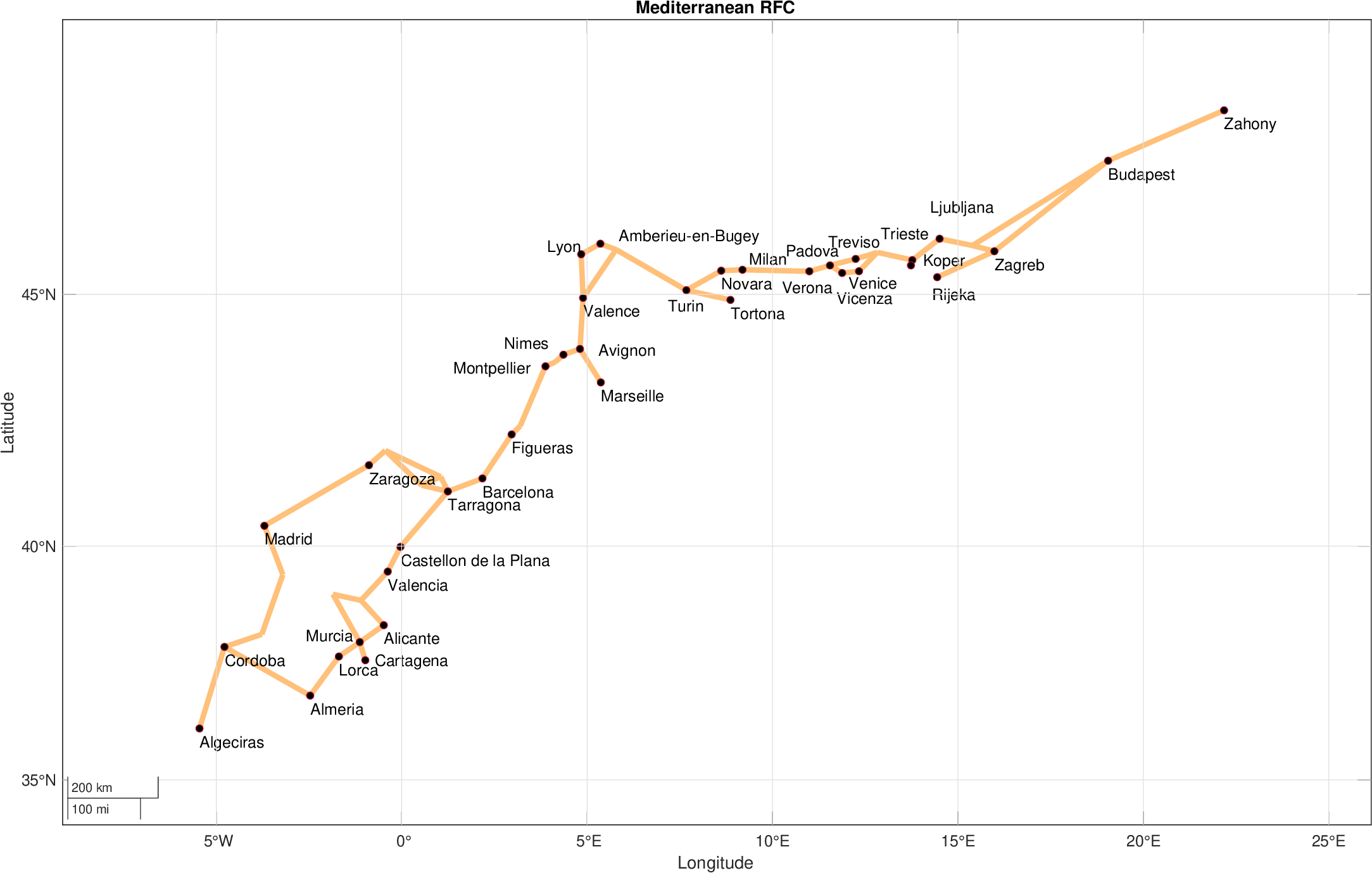}\\
		\includegraphics[width=0.32\textwidth, angle=0, origin=c]{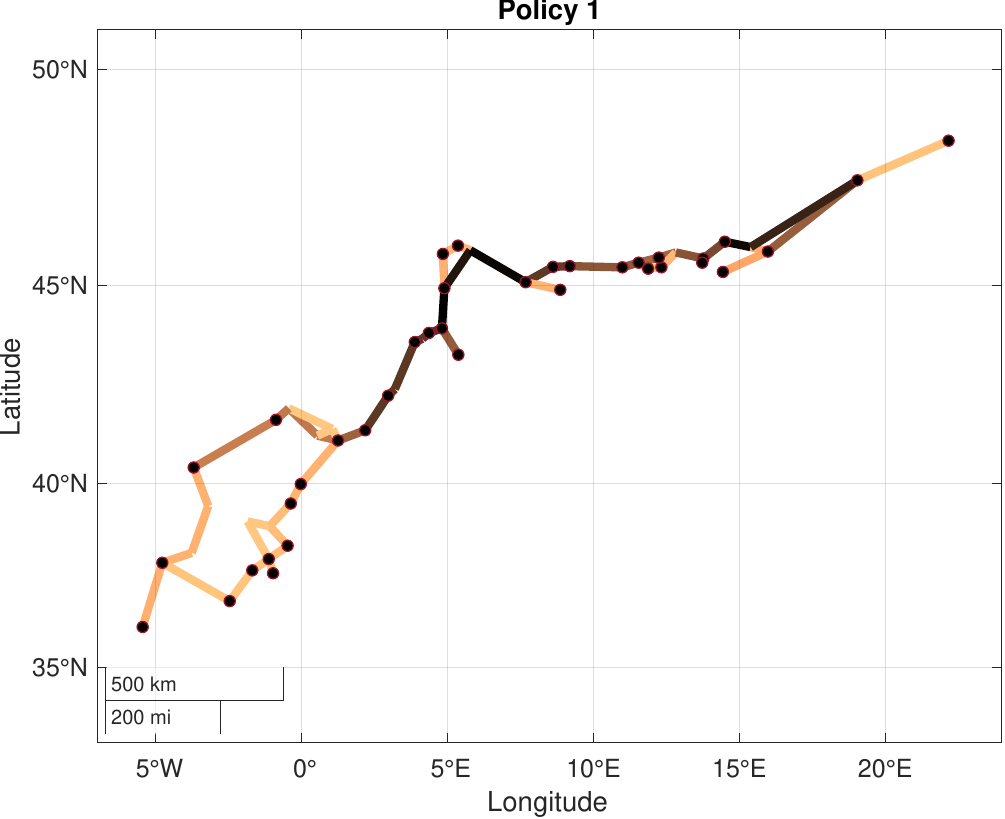}
		\includegraphics[width=0.32\textwidth, angle=0, origin=c]{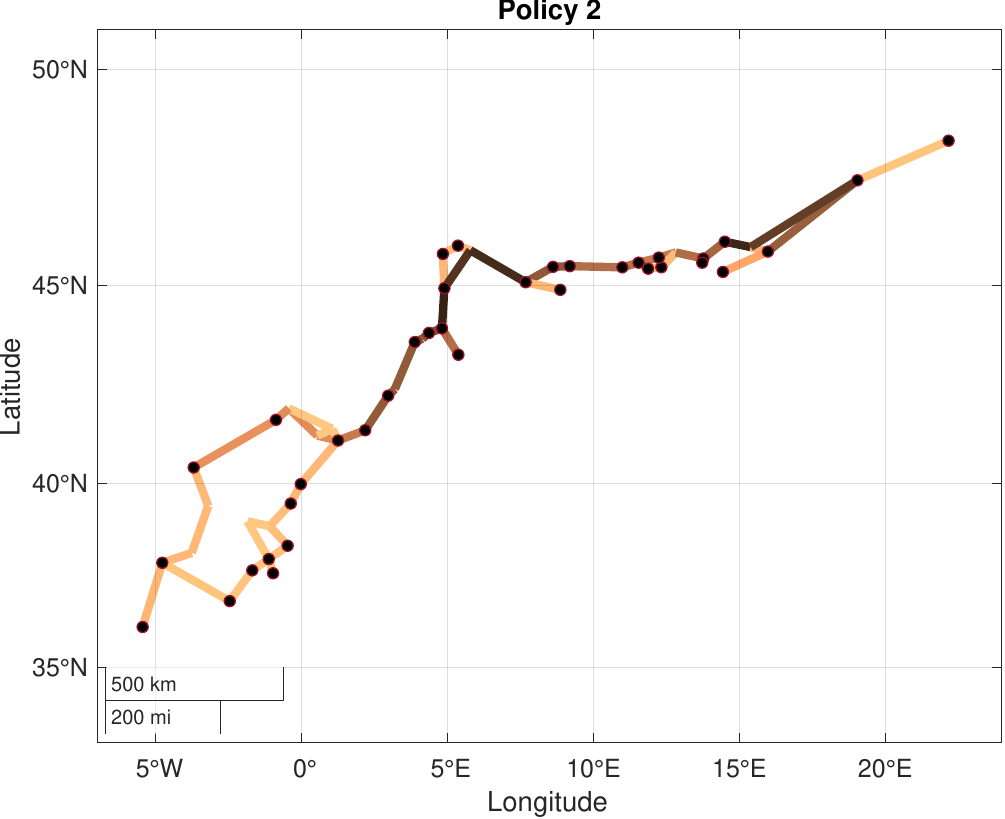}
  \includegraphics[width=0.32\textwidth, angle=0, origin=c]{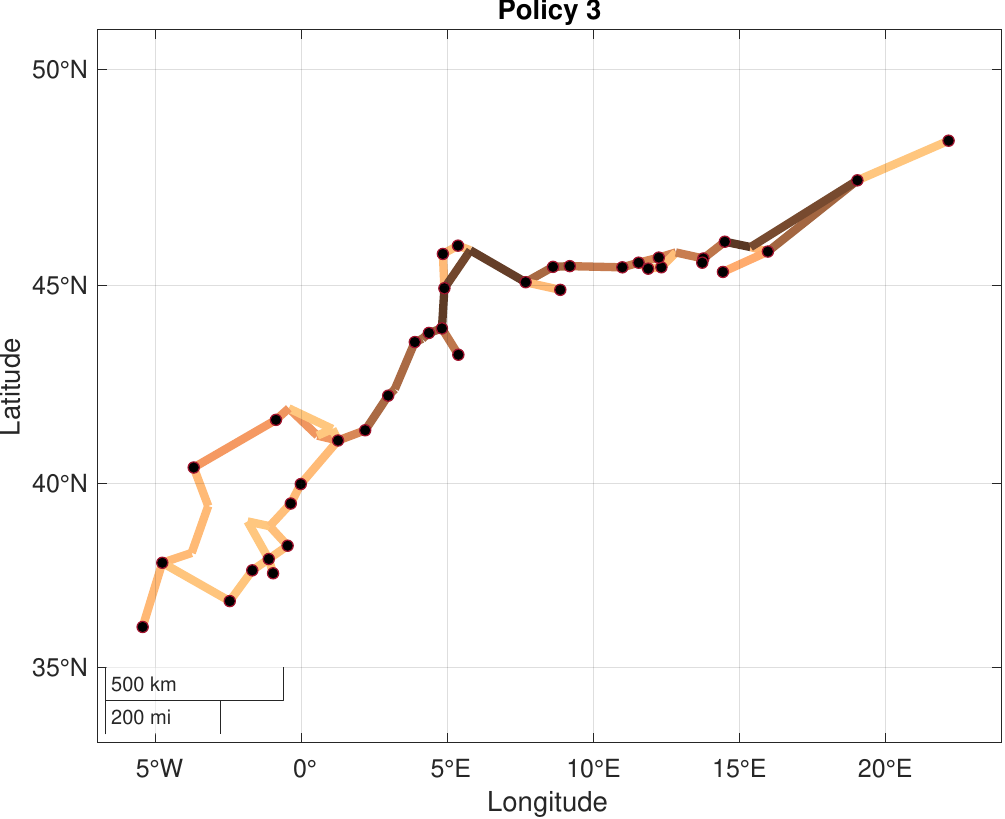}
	\end{center}
 \caption{Optimal flows for the policies \label{fig:flows}}
\end{figure}

The distribution of freight flows is another crucial aspect for infrastructure planning. Figure~\ref{fig:flows} illustrates the optimal freight flows across different policies, showing that traffic distribution is not uniform throughout the network. This insight is crucial for effective capacity management, as it helps identify potential congestion points and underutilized segments.

\begin{figure}[!h]
	\begin{center}
		\includegraphics[width=0.32\textwidth, angle=0, origin=c]{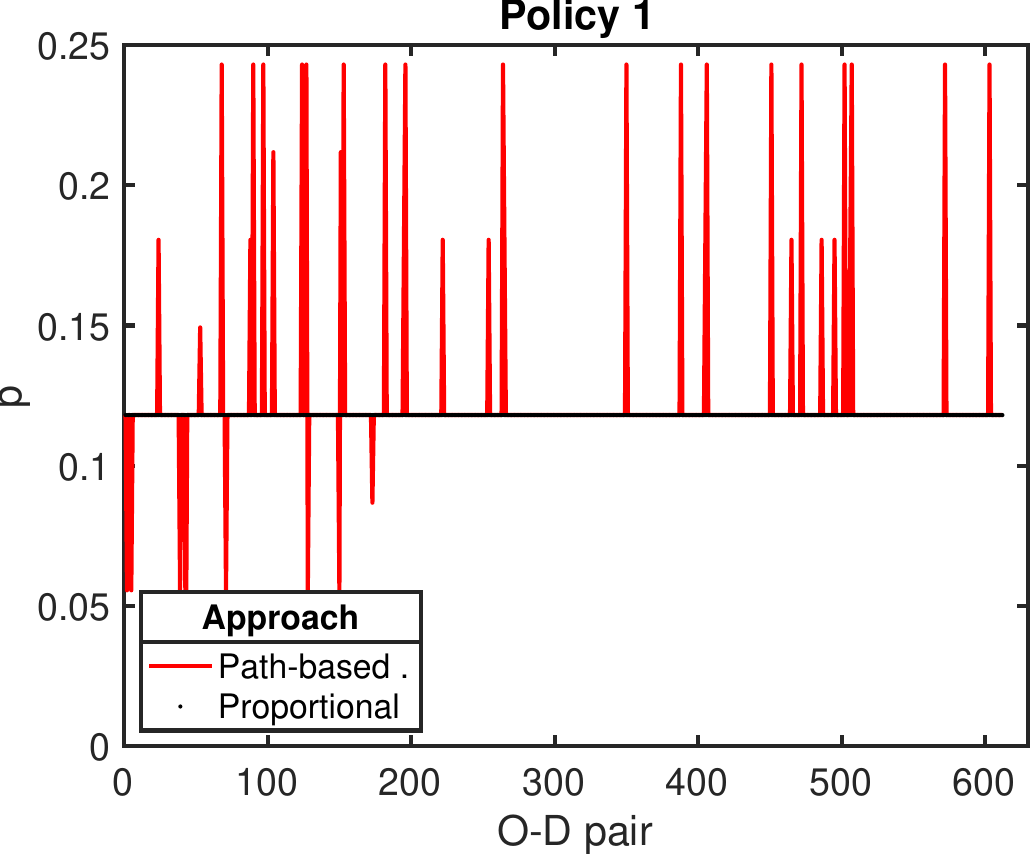}
		\includegraphics[width=0.32\textwidth, angle=0, origin=c]{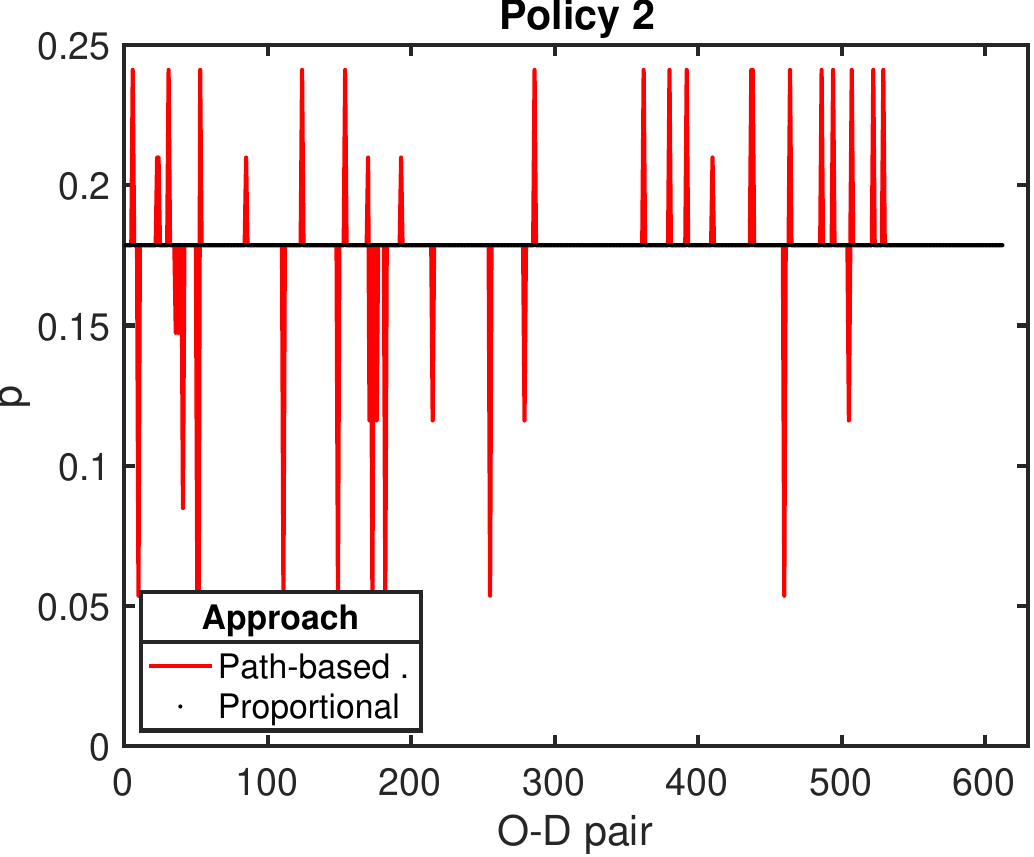}
        \includegraphics[width=0.32\textwidth, angle=0, origin=c]{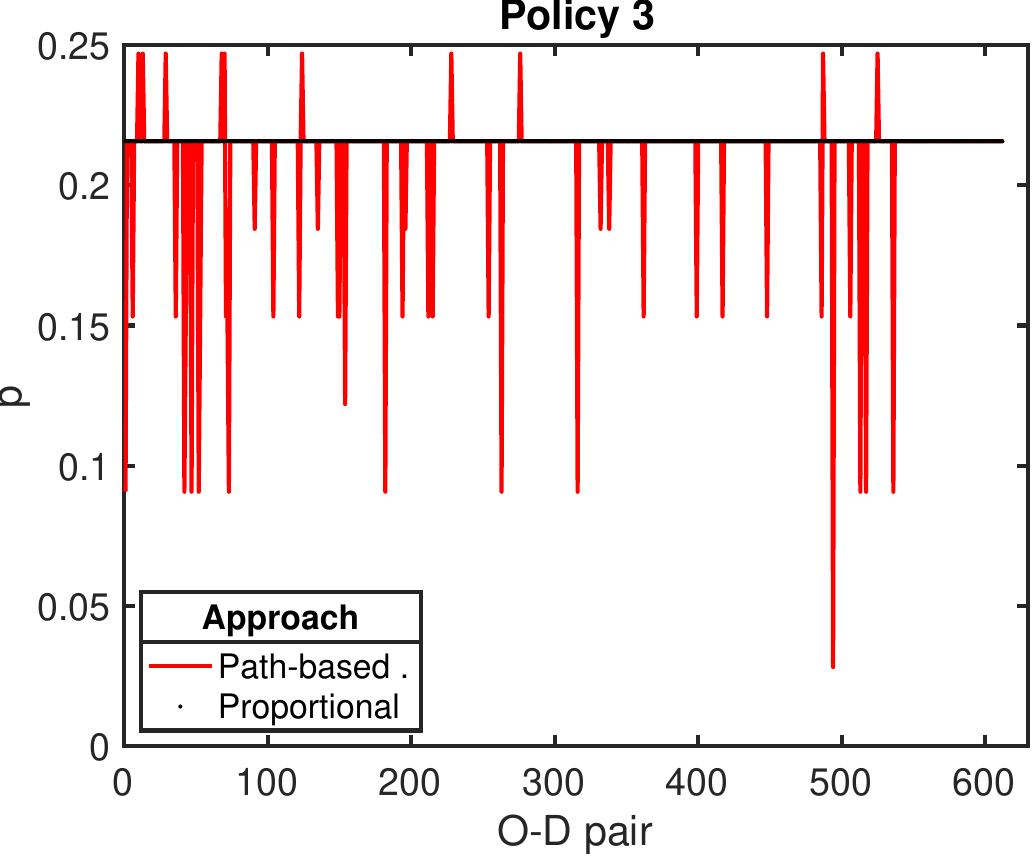}
	\end{center}
 \caption{Solution obtained for the several policies \label{fig:solucion}}
\end{figure}

The impact of pricing strategies is further illustrated in Figure~\ref{fig:solucion}, which shows the computed \glspl{TAC}, revealing two key observations. Firstly, the path-based approach closely aligns with the proportional scheme, with differences confined to a limited subset of origin-destination pairs. Secondly, when environmental costs are given significant weight ({\sl Policy 1}), the applied rates \( p \) amount to approximately $13\%$ of fixed costs—roughly half the rate observed when environmental costs are disregarded, which reaches $22\%$ ({\sl Policy 3}).

In conclusion, the pricing process is inherently multi-objective, as it must balance \gls{IM} revenue generation with broader considerations such as negative externalities and market dynamics. The model provides key performance indicators, summarized in Table~\ref{tab:Results}, which support informed decision making. These indicators include network use levels, revenue generated by the \glspl{FOC}, and the railway's market share, all of which play a role in shaping sustainable and efficient rail policies. While the model allows for comprehensive scenario evaluation, the final decision ultimately rests with the \gls{IM}'s strategic objectives. Specifically, the \gls{IM} must assess whether the revenue reduction observed in {\sl Policy 1} is justified by the associated environmental benefits, increased sector-wide revenues, or enhanced infrastructure use.

\subsection{Experiment II: Establishing track access charges from the state's perspective}

The previous section illustrated the application of the proposed model at a tactical planning level, where the existing railway network was used to estimate \glspl{TAC} for freight transport. Experiment I adopted the perspective of the \gls{IM}, evaluating various pricing schemes and policies based on different metrics. In contrast, Experiment II shifts the perspective to that of the State, the primary investor in railway infrastructure, with the objective of assessing whether transport externalities should be incorporated into \gls{TAC} pricing. This question is central to ongoing economic debates about whether \glspl{TAC} should be designed solely for cost recovery by the \gls{IM} or whether they should be part of a broader public investment strategy. This experiment explores the latter perspective, emphasizing the importance of the integration of externalities in \gls{TAC} determination.

The assessment of railway investment policies and their impacts on modal shift and external freight transport costs is typically treated as a capital investment decision, often analyzed using \gls{CBA}. This methodology compares total project costs over its lifetime with anticipated benefits. However, a full-scale \gls{CBA} for RFC6 is beyond the scope of this paper. Applying the proposed model in a strategic investment context presents two key limitations. Firstly, passenger transport is not explicitly considered, meaning that its indirect benefits are not accounted for. Secondly, the model employs a mesoscopic discrete-event simulation rather than a microscopic one, meaning that it does not compute the expanded effective capacity resulting from infrastructure investments. Instead, the goal of this experiment is to estimate the order of magnitude of costs and benefits, rather than their precise values, to assess the relevance of incorporating externalities into \gls{TAC} pricing.

The baseline year for this analysis is $2016$, when $21.1$ million metric tons (M mt) of freight were transported by rail on RFC6, out of a total of $185.14$ M mt transported across all modes. At that time, the rail market share was $11.41\%$, while road transport dominated with $77.79\%$ (approximately $144.0$ M mt). This study projects conditions for the year $2030$, assuming that the core RFC6 investments have been completed, leading to a fully two-way railway network, commercial train speeds of $100$ km/h, and a line capacity of six trains per hour per direction.

The hypothesis assumes that, without intervention, the rail market share would have remained at $2016$ levels, whereas infrastructure investments lead to an increase. This study evaluates the magnitude of this increase using $2023$ monetary values. From the State's perspective, the following benefits are considered:
\begin{enumerate}
    \item {\sl \glspl{TAC}:} This represents revenue for the \gls{IM}, as determined by the model in Experiment I.
    
    \item {\sl Externalities:} Freight transport generates negative externalities, including air pollution, greenhouse gas emissions, water pollution, noise pollution, accidents, and land use. The study by \cite{DHS15} provides an overview of the external costs associated with these factors. Due to significant variability in influencing parameters, many studies estimate external costs within ranges.
    Table~\ref{tab:externalities} presents the aggregated externality costs for key environmental and societal factors (air pollution, greenhouse gases, noise, accidents, and land use). To adjust these estimates to $2023$ values, a $37.85\%$ increase is applied, corresponding to an annual discount rate of $2.5\%$. Since the external costs of road transport are substantially higher than those of rail, an increase in rail market share results in a reduction of these externalities. For instance, the minimum estimated monetary value of externality savings is computed as:
    \begin{equation}
        1.3785 \times \frac{0.42-0.06}{100} \times \text{tons transferred to rail} \quad (\text{in 2023 €}).
    \end{equation}

    \begin{table}[h]
\centering
\caption{Aggregated cost figures of the negative externalities \label{tab:externalidades}}
\label{tab:externalities}
\resizebox{\textwidth}{!}{
\begin{tabular}{lccp{0.1cm}cc}
\toprule
& \multicolumn{2}{c}{\bf Road}&& \multicolumn{2}{c}{\bf Rail}
\\
\cline{2-3} \cline{5-6} 
\bf Reference & \bf Lower Bound & \bf Upper Bound && \bf  Lower Bound  & \bf Upper Bound\\
  & \bf  \euro ct / \bf ton $\cdot$ km  & \bf \euro ct /  ton $\cdot$ km  && \bf \euro ct / ton $\cdot$ km  & \bf \euro ct / ton $\cdot$ km \\
\hline
\cite{Eco04} & 4.12 & 4.12  && 1.15 & 1.15 \\
\cite{McA10} & 0.05 & 10.95 && 0.05 & 1.49 \\ 
\cite{DeM10} & 0.39 & 20.22 && 0.06 & 0.22 \\ 
\cite{SKS12} & 1.07 & 1.07 && 0.08 & 0.08 \\ 
\cite{VTPI}& 0.17 & 3.03 && 0.04 & 1.17 \\ 
\hline
\bf Average (this study)& 0.42 & 8.82 && 0.06 & 0.74 \\ 
\bottomrule
\end{tabular}
}
\end{table}


    \item {\sl \gls{FOC} Benefits:} Infrastructure improvements lead to higher commercial speeds, which in turn reduce operational costs for \glspl{FOC}. These cost savings are considered a direct benefit for railway freight operators.
    
    \item {\sl Social Benefit:} Rail transport costs are generally lower than road transport costs. Transferring freight from road to rail reduces overall transportation costs, benefiting society as a whole. However, this modal shift also has economic implications, particularly in tax revenue. Since road transport generates more tax revenue per ton-km than rail transport, a reduction in road freight activity affects the State's tax income. To account for this effect, a tax revenue reduction factor of $32\%$ is applied to both the social benefit estimation and effective State investment. This reflects the portion of the State's investment that is indirectly recovered through tax revenue from the resulting economic activity.
\end{enumerate}

Table~\ref{tab:experiment_2_benefits} provides estimates for the various benefits derived from the implementation of \glspl{TAC} and the inclusion of externalities. Since externalities are assessed within a range, both minimum and maximum values are given. As a result, the total benefit, calculated as the sum of all components—including externalities—is also expressed as an interval.

A key observation is that total benefits are significantly higher when externalities are integrated into the \gls{TAC} calculation compared to policies where they are not considered. This difference is primarily driven by externality reductions and social-cost adjustments, which reflect the effects of modal shift from road to rail. These components play a crucial role in offsetting the reduction in \gls{TAC} revenue, highlighting the importance of incorporating external costs into railway pricing strategies.

The operational profit for the \glspl{FOC} accounts for the efficiency gains from improved railway operations, specifically attributed to increased commercial speed. However, as railway operations expand, the network experiences higher congestion levels, leading to a potential reduction in commercial speed. The results in Table~\ref{tab:experiment_2_benefits} suggest that these two opposing effects—increased traffic and congestion-related speed reduction—largely counterbalance each other. Consequently, the impact of \gls{TAC} structures on commercial speed does not appear to be a decisive factor in determining the most effective pricing scheme.

\begin{table}[h]
\centering
\caption{Benefits (in M \euro) vs scenarios\label{tab:experiment_2_benefits}}
\resizebox{\textwidth}{!}{
\begin{tabular}{llp{1.5cm}p{0.75cm}p{0.75cm}p{1.5cm}p{1.5cm}p{0.75cm}p{0.75cm}}
\toprule
&& \bf  \glspl{TAC}  & \multicolumn{2}{l}{\bf Externalities} & \bf  \gls{FOC}'s benefit  & \bf Social's benefit  & \multicolumn{2}{c}{\bf Total benefit} \\ 
\toprule
\multirow{ 2}{*}{\bf Policy 1}& Path-based & 105.51 & [91.00, & 730.07] & 88.33 & 286.33 & [571.16, & 1210.23] \\ 
& Proportional & 101.69 & [90.16, & 723.31] & 88.9 & 292.45 & [573.2, & 1206.35] \\ 
\hline 
\multirow{ 2}{*}{\bf Policy 2}& Path-based &  123.39 & [47.18, & 378.54] & 88.5 & 132.72 & [391.79, & 723.15] \\ 
&Proportional& 121.4 & [43.14, & 346.12] & 85.75 & 118.57 & [368.87, & 671.84] \\ 
\hline 
\multirow{ 2}{*}{\bf Policy 3}&Path-based &128.38 & [22.23, & 178.36] & 84.85 & 60.92 & [296.38, & 452.51] \\ 
&Proportional &125.27 & [16.49, & 132.32] & 80.8 & 44.37 & [266.93, & 382.77] \\ 
\bottomrule
\end{tabular}
}
\end{table}

\begin{figure}[h]
    \centering
    \begin{minipage}{0.45\textwidth}
        \centering
        \resizebox{0.65\textwidth}{!}{ 
            \begin{tabular}{ll}
                \toprule
                \bf Period & \bf \euro (billion)\\
                \toprule
                2014-2016& 1.362\\
                2017-2020& 8.523\\
                2021-2025 & 30.447\\
                2026-2030 & 72.415\\
                after 2030 & 32.498\\
                unknown  & 1.123\\
                \hline
                \bf Total &  146.368\\
                \bottomrule
            \end{tabular}
        }
       
    \end{minipage}
    \hfill
    \begin{minipage}{0.45\textwidth}
        \centering
        \includegraphics[width=\textwidth]{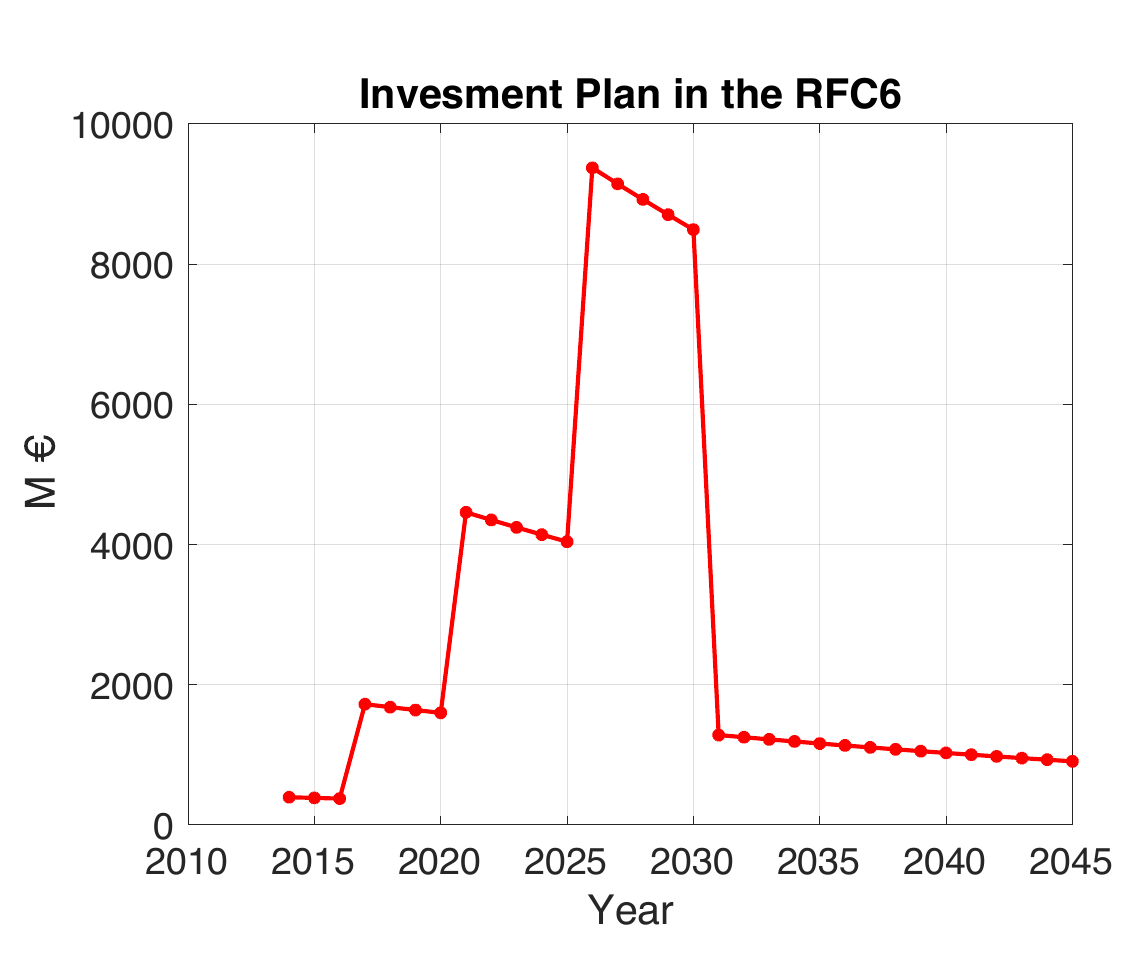}
    \end{minipage}
    \caption{Investment plan in the RFC6 (source \cite{Rad22}) and its net present value (NPV)\label{fig:budget}}
\end{figure}

Figure~\ref{fig:budget} shows the historical and planned investments for RFC6, with each value representing the amount allocated for the corresponding year. The \gls{NPV} methodology is applied to express the current (as of $2023$) value of future cash flows from these infrastructure investments. To compute the \gls{NPV}, the expected timing and magnitude of future cash flows are estimated, and a discount rate of $2.5\%$ is applied, representing the minimum acceptable rate of return for these public investments. Figure~\ref{fig:budget} illustrates the investment plan in \gls{NPV} terms, under the assumption that the State will recover $32\%$ of the investment through tax revenues. Based on this projection, the average annual investment required for RFC6 amounts to $2880.4$ M\euro\, (expressed in $2023$ euros).

Assuming an annual investment cost of $2880.4$ M\euro, the Benefit–Cost Ratio (BCR) is calculated to assess the relative scale of investment costs to expected benefits. It is important to note that this analysis does not represent a comprehensive \gls{CBA}, as several critical factors have been omitted, including the residual value of investments, the temporal distribution of costs and benefits, and the exclusion of passenger transport impacts, among others. The purpose of the BCR in this context is to provide an approximate measure for comparing different policies and evaluating the role of \glspl{TAC} in shaping investment decisions.

Table~\ref{tab:experiment_2_CBA} presents the estimated BCR indices. The inclusion of externalities in \gls{TAC} policy leads to a $38\%$ to $147\%$ increase in the lower bound of the BCR and a $75\%$ to $210\%$ increase in the upper bound. These results highlight a key policy trade-off: if \glspl{TAC} are viewed purely as a cost-recovery mechanism, then {\sl Policy $3$}, which maximizes \gls{IM} revenue, would be the optimal choice. However, if \glspl{TAC} are considered within a broader investment framework, {\sl Policies} $1$ and $2$, which account for environmental externalities, would be preferable.

\begin{table}[h]
\centering
\caption{Benefit Cost Ratio\label{tab:experiment_2_CBA}}
\resizebox{0.8\textwidth}{!}{
\begin{tabular}{llcc}
\toprule
&& \bf Lower Bound ($\%$)   & \bf Upper Bound ($\%$) \\
\toprule
\multirow{ 2}{*}{\bf Policy 1}& Path-based &  19.83 & 42.02\\
& Proportional & 19.9 & 41.88 \\ 
\hline 
\multirow{ 2}{*}{\bf Policy 2}& Path-based & 13.6 & 25.11 \\ 
&Proportional & 12.81 & 23.32 \\ 
\hline 
\multirow{ 2}{*}{\bf Policy 3}&Path-based & 10.29 & 15.71 \\ 
&Proportional &9.27 & 13.29 \\ 
\bottomrule
\end{tabular}
}
\end{table}


While the BCR index is not fully estimated in this study, it provides valuable insight into public investment strategies. In private sector investments, a BCR exceeding $100\%$ is typically required to justify funding. However, for public investments, this threshold is not necessarily applicable, as social, environmental, and economic policy objectives must also be considered. For instance, the EU's White Paper ``Roadmap to a Single European Transport Area: A Competitive and Resource Efficient Transport System'' establishes a target that $30\%$ of long-distance freight transport (over $300$ km) currently handled by road should be shifted to rail or waterborne transport by $2030$. In the case of RFC6, achieving this $30\%$ modal shift would require rail and maritime transport to absorb $23.33\%$ of the total freight market by $2030$. The analysis suggests that corridor improvements alone are expected to achieve a $6\%$ modal shift, indicating that additional measures will be required to meet EU environmental and transport targets.

A comprehensive policy package will likely be necessary, combining measures such as reducing port and rail fees and increasing road transport taxes to further encourage modal shift. The BCR framework can serve as a valuable tool for prioritizing these interventions, helping to determine which policy measures should be implemented to ensure the EU's transport and sustainability objectives are met.

Transport infrastructure investment and pricing policies must balance financial sustainability, environmental objectives, and market competitiveness. Achieving sustainable transport requires a combination of measures, much like the energy sector relies on diverse sources to meet demand at different costs. In this context, rail track access charging systems should incorporate externalities and set prices that, while potentially higher than those of other transport policies, are crucial for facilitating the necessary modal shift.

\section{Conclusions}

This study introduces a novel methodology for pricing railway infrastructure access in deregulated railway markets. The proposed approach is based on a dynamic freight flow model that integrates a logit-based modal split function, incorporating \gls{TAC} and network congestion effects (i.e., travel time). Upon discretization, the model transforms into a mesoscopic simulation framework, enabling railway traffic simulations that account for capacity constraints. This approach facilitates the development of dynamic pricing strategies that optimize \gls{TAC} revenue for the \gls{IM} while incorporating the economic value of reducing negative externalities.

The methodology was applied to the Mediterranean Rail Freight Corridor, highlighting the computational challenges associated with optimizing the simulation model. Due to the discrete-event nature of the simulation, gradient-based optimization methods proved ineffective, necessitating the use of heuristic algorithms such as \gls{PS}. For the path-based approach, the \gls{PS} algorithm required over $10$ hours of computational time to converge to an optimal solution.

Additionally, the so-called proportional approach was analyzed, wherein a single proportionality constant was applied to price all train paths. This alternative method required less than a minute of computation time, offering a transparent and easily implementable solution for \glspl{FOC}, albeit with a slight reduction in the optimality of the results obtained.

The findings from Experiment I indicate that pricing strategies incorporating externalities lead to a moderate reduction in the \gls{IM}'s direct economic return, but generate significant environmental benefits. Experiment II further examined whether \gls{TAC} should be treated as a cost-recovery mechanism or as an investment tool. The results suggest that when \gls{TAC} is framed as an investment problem, the economic benefits for the \gls{IM} are significantly lower than the social and environmental externalities. This underscores the necessity of integrating externalities into railway pricing policies. Overall, the study suggests that well-designed \gls{TAC} policies can support the transition toward a more sustainable and energy-efficient transport system, aligning with EU climate goals and transport sector emission reduction targets.

The proposed model operates at a tactical planning level and incorporates certain simplifications regarding railway flow dynamics. Future research should focus on enhancing the model's granularity by integrating microscopic simulation approaches to improve railway system capacity assessments. Furthermore, expanding the model to explicitly account for passenger transport and multimodal interactions could provide a more comprehensive analysis of network performance and the broader implications of \gls{TAC} pricing strategies.

\section*{CRediT authorship contribution statement}
\textbf{Ricardo García-Ródenas:} Conceptualization, Methodology, Software, Data curation, Writing – original draft, Supervision.
\textbf{Esteve Codina}: Conceptualization, Methodology, Writing – original draft.
\textbf{Luis Cadarso}: Conceptualization, Methodology, Writing – original draft.
\textbf{María Luz López-García}: Conceptualization, Methodology, Writing – review \& editing.
\textbf{José Ángel Martín-Baos}: Conceptualization, Software, Data curation, Writing – review \& editing.


\section*{Acknowledgements}
This work was supported by grants PID2020-112967GB-C31, PID2020-112967GB-C32,  PID2020-112967GB-C33, and TED2021-130347B-I00 funded by MCIN/AEI/10.13039/501100011033 and by {\sl ERDF A way of making Europe}. Additional funding was provided by grant 2022-GRIN-34249 from the University of Castilla-La Mancha and ERDF.

\bibliographystyle{elsarticle-harv} 
\bibliography{references}

\newpage

\end{document}